# All Finite (Anti)Hermitian Irreducible Representations of the de Sitter and Anti-de Sitter Lie Algebras and Their Lorentz Structure

R.A.W. Bradford

University of Bristol, Faculty of Engineering

## Abstract

Because of the importance of unitarity in quantum physics, work on the representations of the de Sitter group has focussed on the unitary case, which necessarily means infinite dimensional matrices for this non-compact group. Here we address the finite dimensional representations resulting from the requirement that the Lie algebra generators are either Hermitian or anti-Hermitian. The complete classification of all such irreducible representations is found and their matrix elements specified. These irreducible representations (irreps) are based on "backbones" defined as the homogeneous Lorentz sub-algebra and consisting of direct sums of the finite irreps of the homogeneous Lorentz algebra (HLA). All these irreps have been characterised by a diagrammatic depiction conveniently labelled by two positive integers $m/n$ where $m, n \in \{1,2,3, \dots\}$ excepting only $m = n = 1$. The two Casimir invariants have been given explicitly in terms of the two positive integers, $m, n$, equs.(11.11a,b), and shown to reproduce the result familiar from the study of the infinite dimensional, unitary representations $-C1 = p(p+1) + (q+1)(q-2)$ and $-C2 = p(p+1)q(q-1)$ but with a reduced set of possible $p, q$ values. There are 17 irreps with dimensions less than 105, having dimensions 4, 5, 10, 14, 16, 20, 30, 35(x2), 40, 55, 56, 64, 80, 81, 84 and 91. The corresponding irreps of the anti-de Sitter group follow immediately from the replacement of the 4-momentum operators thus $V_\mu \to iV_\mu$.

## Contents

## 1. Purpose

The representations of the de Sitter group (dSG) and the anti-de Sitter group (AdSG) have received considerable attention in the past, initially by Thomas, Ref.[1], which was subject to some later corrections by Newton, Ref.[2], and consolidated by Dixmier, Ref.[3]. Most recently Enayati et al have presented the topic in comprehensive book form, Ref.[4]. However, because of the importance of unitarity in quantum mechanics, the focus has been on unitary representations. But the dSG and AdSG are non-compact and hence admit only infinite dimensional unitary representations. Consequently, existing work has focussed on the infinite dimensional unitary representations. From the perspective of the Lie algebra, unitarity of the group arises from all the generators of the Lie algebra being Hermitian.

In this paper we present the complete set of irreducible finite dimensional representations of the dSA and AdSA subject to the less restrictive condition that the generators are either Hermitian (H) or anti-Hermitian (AH). (Note that with this requirement to be H or AH, irreducibility and indecomposability are the same). Expressions are provided for the matrix components of all the generators.

The homogeneous Lorentz algebra (HLA) is a sub-algebra of the dSA and the AdSA. We are interested in the Lorentz structure of the irreducible finite representations of the dSA and AdSA; in other words, to what set of representations of the HLA does an irreducible representation (irrep) of dSA/AdSA reduce? This structure is automatically provided if the finite irreps of the dSA/AdSA are built from a "backbone" consisting of direct sums of the known finite irreps of the HLA (see (4.5), below, for details).

Existing work on the infinite dimensional unitary representations of the dSG has allowed all such representations to be parametrised in terms of an integral, or half-integral, parameter $p$, and another parameter, $q$. The latter may take integral or half-integral values, or, for some classes of representation, values from a continuous range. The two Casimir invariants, which characterise the irreps, are given by simple expressions in terms of $p$ and $q$. The same expressions are found to apply also in the finite dimensional case, restricted to H or AH generators, but occur only for the cases $q = 0$ or $q = p$. Thus, the complete set of all irreps in the finite, H or AH, case is severely restricted compared with the set of irreps available in the infinite dimensional unitary case – even if attention is confined to the discrete series of the latter.

Throughout we use the term "half integral" to mean an odd integer / 2 (hence, if referring to spin, "half-integral" would mean fermionic). To avoid repeating that the finite irreps sought here are subject to the condition that all generators are Hermitian or anti-Hermitian, please take this as understood throughout.

## 2. The Lie Algebra

It will be convenient to express the dSA/AdSA algebra in terms of the familiar rotation generators, $J_i$, the generators of boosts, $K_i$, and the four generators of displacements, $V_\mu$. Latin subscripts take the values $x, y, z$ whilst Greek subscripts take values $t, x, y, z$. The $J_i, K_i, V_\mu$ generate the group of isometries of a 4-sphere in (4,1) space in the dSA case, or in (3,2) space in the AdSA case, i.e., group elements are transformations that preserve the 5-dimensional metric $t^2 \mp w^2 - (x^2 + y^2 + z^2)$ where the upper sign refers to the dSA and the lower sign to the AdSA.

The HLA is a sub-algebra whose commutation relations (CRs) are,

$$[J_i, J_j] = i\epsilon_{ijk}J_k \qquad [K_i, K_j] = -i\epsilon_{ijk}J_k \qquad [J_i, K_j] = i\epsilon_{ijk}K_k \tag{2.1}$$

The remaining CRs all involve the $V_\mu$, the generalisation of the usual 4-momentum operators. They may be derived from the stated isometry condition. Importantly, the CRs between $V_\mu$ and $J_i$ or $K_i$ are equivalent to requiring that the four quantities $V_\mu$ transform as a Lorentz 4-vector (thus justifying the notation). Alternatively, the CRs imbued by the isometry condition imply that $V_\mu$ is a Lorentz 4-vector. These CRs are,

$$[J_i, V_j] = i\epsilon_{ijk}V_k \qquad [K_i, V_j] = -iV_t\delta_{ij} \qquad [J_i, V_t] = 0 \qquad [K_i, V_t] = -iV_i \tag{2.2}$$

So far, (2.1) and (2.2) are identical to the equivalents in the Poincaré Lie algebra. However, in flat Minkowski space the 4-momentum operators, i.e., the generators of displacements, are mutually commuting. The hallmark of intrinsic curvature is that the generators of displacements do not commute. For the dSA/AdSA geometries we find,

$$[V_i, V_j] = \pm i\epsilon_{ijk}J_k \qquad [V_t, V_i] = \pm iK_i \qquad (+ = \text{dS},\ - = \text{AdS}) \tag{2.3}$$

The first CR in (2.1) is not consistent with the $J_i$ being anti-Hermitian (AH). So the $J_i$ must be Hermitian (H) – for our purposes, i.e., because we are interested only in the case that all generators are H or AH.

All the CRs in (2.1), (2.2) and (2.3) are consistent with all the other generators, $K_i$ and $V_\mu$, being H or all AH (or certain mixtures). However, no finite representations will exist with all generators H; some must be AH. This is familiar from the well-known classification of all the finite representations of the HLA sub-algebra, i.e., (2.1). The corresponding group is also non-compact so the finite representations must involve some generators which are AH. But the $J_i$ must be H, so the $K_i$ must be AH. This is explicitly manifest by considering the complexification of the algebra and defining,

$$A_i = (J_i + iK_i)/2 \qquad B_i = (J_i - iK_i)/2 \tag{2.4a}$$

i.e.,

$$J_i = A_i + B_i \qquad K_i = i(B_i - A_i) \tag{2.4b}$$

In terms of which the CRs (2.1) can be written in decoupled form as,

$$[A_i, A_j] = i\epsilon_{ijk}A_k \qquad [B_i, B_j] = i\epsilon_{ijk}B_k \qquad [A_i, B_j] = 0 \tag{2.5}$$

In other words, the HLA is isomorphic to the direct sum of two copies of the (complexified) $\mathfrak{su}(2)$. But CRs (2.5) imply that both of $A_i$ and $B_i$ must be H, so that (2.4b) implies that the $K_i$ must be AH, as we have already concluded.

It can thus be seen that the CRs in (2.1), (2.2) and (2.3) are consistent with $V_i$ being H if $V_t$ is AH, or vice-versa. These cases align with dS and AdS respectively. This can be seen by noting that, as well as the 3-rotations, $\{J_i\}$ and the HLA, $\{J_i, K_i\}$ being sub-algebras, the $\{J_i, V_i\}$ also provide a sub-algebra. If we choose the dSA sign in (2.3), i.e., $[V_i, V_j] = i\epsilon_{ijk}J_k$, then $\{J_i, V_i\}$ generates SO(4), i.e., rotations in $E^4$, which requires the $V_i$ to be H to generate these finite unitary transformations. Conversely, if we choose the AdS sign in (2.3), i.e., $[V_i, V_j] = -i\epsilon_{ijk}J_k$, then the $\{J_i, V_i\}$ provide an algebra identical to the HLA with $K_i$ replaced by $V_i$ and hence requiring the $K_i$ to be AH.

Note that in both cases the $\{J_i, V_i\}$ provide a sub-algebra which is either isomorphic to $\mathfrak{su}(2) \oplus \mathfrak{su}(2)$ or its complexification is isomorphic to $\mathfrak{su}(2) \oplus \mathfrak{su}(2)$. This has formed the basis of some approaches to the representations of the dSA, using as the "backbone" direct sums of the known irreps of $\mathfrak{su}(2) \oplus \mathfrak{su}(2)$, i.e., the SO(4) subgroup. This alternative means of building the irreps will be important in a later part of this paper (§11).

To summarise we have,

- $A_i$ and $B_i$ are Hermitian.
- $J_i$ are Hermitian
- $K_i$ are Anti-Hermitian.
- dS: $V_i$ are Hermitian, $V_t$ is Anti-Hermitian.
- AdS: $V_i$ are Anti-Hermitian, $V_t$ is Hermitian.

## 3. The Anti-de Sitter Case

All irreps of the AdSA can be formed from those of the dSA by the replacement $V_j \to iV_j$ and $V_t \to iV_t$,noting that this leaves all CRs (2.1) and (2.2) invariant, whereas the RHS of (2.3) change sign, as required. We therefore confine attention to the dS case, the AdS case being understood.

## 4. The Irreps of the HLA

These are well known, following from those of $\mathfrak{su}(2)$, and we simply state the result without elaboration. Exactly one irrep of $\mathfrak{su}(2)$ exists, up to equivalence, for every positive integral dimension, 1, 2, 3…. We write this dimension as $(2A+1)$ where $A$ takes all positive indefinite integral or half-integral values. (Throughout we use the term "half integral" to mean an odd integer / 2). Consider the following combinations,

$$J^{\pm} = J_x \pm iJ_y \qquad K^{\pm} = K_x \pm iK_y \qquad A^{\pm} = A_x \pm iA_y \qquad B^{\pm} = B_x \pm iB_y \tag{4.1}$$

Indexing the rows/columns of the representation by indices $a \in \{-A, -A+1, \dots. A-1, A\}$ the representations are,

$$\left(A^{+(A)}\right)_{a_1a_2} = r_{a_2}^{(A)}\delta_{a_1a_2+1} \tag{4.2a}$$

$$\left(A^{-(A)}\right)_{a_1a_2} = s_{a_2}^{(A)}\delta_{a_1a_2-1} \tag{4.2b}$$

$$\left(A_z^{(A)}\right)_{a_1a_2} = a_1\delta_{a_1a_2} \tag{4.2c}$$

where,

$$r_{a_2}^{(A)} = \sqrt{(A-a_2)(A+a_2+1)} \tag{4.2d}$$

$$s_{a_2}^{(A)} = \sqrt{(A+a_2)(A-a_2+1)} \tag{4.2e}$$

The same expressions apply, of course, to $B_i$. Note that although the indices $a_1$, etc., may be integral or half-integral, adjacent indices differ by 1.

When the abstract generators are realised as matrices, the additions in (2.4a,b) must be understood in terms of the addition of implied direct products, that is, for example,

$$J_i = A_i \otimes \mathbb{I}_{2B+1} + \mathbb{I}_{2A+1} \otimes B_i \tag{4.3}$$

where $\mathbb{I}_{2A+1}$ and $\mathbb{I}_{2B+1}$ are unit matrices of dimension equal to that of the $A_i$ and $B_i$ matrices respectively. It is convenient when dealing with direct products of matrices to use a double index notation. Hence the $J_i$ and $K_i$ are represented by the following matrices,

$$\left(J^{+(A,B)}\right)_{a_1 b_1, a_2 b_2} = r_{a_2}^{(A)} \delta_{a_1 a_2+1} \delta_{b_1 b_2} + r_{b_2}^{(B)} \delta_{a_1 a_2} \delta_{b_1 b_2+1} \quad (4.4a)$$

$$\left(J^{-(A,B)}\right)_{a_1 b_1, a_2 b_2} = s_{a_2}^{(A)} \delta_{a_1 a_2-1} \delta_{b_1 b_2} + s_{b_2}^{(B)} \delta_{a_1 a_2} \delta_{b_1 b_2-1} \quad (4.4b)$$

$$i\left(K^{+(A,B)}\right)_{a_1 b_1, a_2 b_2} = r_{a_2}^{(A)} \delta_{a_1 a_2+1} \delta_{b_1 b_2} - r_{b_2}^{(B)} \delta_{a_1 a_2} \delta_{b_1 b_2+1} \quad (4.4c)$$

$$i\left(K^{-(A,B)}\right)_{a_1 b_1, a_2 b_2} = s_{a_2}^{(A)} \delta_{a_1 a_2-1} \delta_{b_1 b_2} - s_{b_2}^{(B)} \delta_{a_1 a_2} \delta_{b_1 b_2-1} \quad (4.4d)$$

$$\left(J_z^{(A,B)}\right)_{a_1 b_1, a_2 b_2} = (a_1 + b_1) \delta_{a_1 a_2} \delta_{b_1 b_2} \quad (4.4e)$$

$$i\left(K_z^{(A,B)}\right)_{a_1 b_1, a_2 b_2} = (a_1 - b_1) \delta_{a_1 a_2} \delta_{b_1 b_2} \quad (4.4f)$$

Hence, the irreps of the HLA are uniquely labelled by the pair, $A, B.$ , each of which is an integer or half integer. Every finite irrep of the dSA includes a representation of the HLA which must be equivalent to a finite direct sum of HLA irreps which can be denoted,

$$(A_1, B_1) \oplus (A_2, B_2) \oplus (A_3, B_3) \oplus \ldots \oplus (A_N, B_N) \quad (4.5)$$

for some finite positive integer, $N$, which is the number of HLA irreps involved in the dSA irrep. $N$ will also be referred to as the number of "blocks" because in matrix notation each of the HLA irreps, $(A_i, B_i)$, is a block in the block-diagonal matrix representation defined by (4.5) with (4.4a-f).

The direct sum, (4.5), will be called the "backbone" of the irrep of the dSA. The order of the terms (blocks) in the direct sum makes no difference, of course. Also reversing $A_i \leftrightarrow B_i$ for all $i$ merely changes the sign of $K_i$ and an equivalent representation then results if the sign of $V_t$ is also reversed (all CRs 2.1-3 being invariant under this transformation).

## 5. Structure of the Classification of the Irreps

This paper presents a classification scheme for all the finite irreps of the dSA/AdSA. The classification scheme will be based on specifying the backbones, (4.5), which provide irreps. This has the benefit of displaying the homogeneous Lorentz structure explicitly. The full set of representation matrices will be constructed explicitly for some of the smallest dimension irreps. As defined by (4.5) and (4.4a-f), these representations will be in a basis in which the $J_i$ and $K_i$ are in block-diagonal form. Therefore the representations of $V_\mu$ must involve some off-block-diagonal terms if the representation is to be an irrep.

The development in outline will be,

1) Proof that there is no representation of the dSA based on a backbone with a single block. A proof has been presented by Shurleff in Ref.[5] but will be reiterated in §6, below, for completeness.

2) Proof that non-zero off-block-diagonal (off-BD) terms in (the representations of) $V_\mu$, say the off-BD block between blocks $P$ and $Q$, can occur only if $|A_P - A_Q| = \frac{1}{2}$ and $|B_P - B_Q| = \frac{1}{2}$. Blocks $P$ and $Q$ would then be said to be "compatible". This has also been proved in Ref.[5] but will be reiterated here for completeness, a simple extension of Ref.[5] being to show that the result applies to an arbitrary number of blocks in the backbone. Ref.[5] also derived explicit expressions for the off-BD terms which we will summarise and use extensively in the rest of this paper. These issues are dealt with in Appendix A.

Definition: "Connected"

We have already defined "compatibility" between two blocks $P$ and $Q$, i.e., that $A_P = A_Q \pm \frac{1}{2}$ and $B_P = B_Q \pm \frac{1}{2}$ (signs unrelated). This is a necessary condition for the off-BD block $PQ$ in the representation of the $V_\mu$ to be non-zero. However, compatibility does not mean that the off-BD block $PQ$ in the representation of the $V_\mu$ must be non-zero. If the off-BD blocks $V_{\mu PQ}$ and $V_{\mu QP}$ are actually non-zero, the blocks $P$ and $Q$ are said to be "connected".

Definition: "Monotonic"

An ordered sequence of three or more backbone blocks is said to be "monotonic" if both the $A$ and the B parameters vary monotonically over the sequence. Hence, referring to the signs in $A_P = A_Q \pm \frac{1}{2}$ and $B_P = B_Q \pm \frac{1}{2}$, the possible monotonic sequences of three blocks are $++, ++$ and $--, --$ and $+-, +-$ and $-+, -+$.

Definition: Unique Non-Monotonic 3-Sequence

A sequence of connected blocks from block $I$ to block $J$ is a non-monotonic sequence if it is connected but not monotonic. A 3-sequence involves just three blocks, $IKJ$. A non-monotonic 3-sequence is a unique non-monotonic 3-sequence if there is no other non-monotonic 3-sequence connecting block $I$ to block $J$.

3) Proof that the smallest $A_i$ and the smallest $B_j$ must be zero for the backbone to provide an irrep. This is addressed in §7.
4) Proof that if a proposed backbone-plus-connectivity includes a unique non-monotonic 3-sequence then it cannot provide a representation. This is addressed in §8.
5) The diagrammatic classification of all the irreps of the dSA/AdSA (§9).
6) Demonstration that if a proposed backbone includes a duplicate block, i.e., two or more instances of a given $(A, B)$, then it is not an irrep (either not a representation or is reducible). This is addressed in Appendix C.
7) Explicit construction of all the irreps based on a fully monotonic backbone (involving two distinct infinite classes). This is addressed in Appendix B.
8) Explicit construction of a few of the smallest dimensional irreps based on backbones which are not fully monotonic. This is addressed in §10.
9) Derivation of the Casimir invariants for each of the irreps thus classified. This is addressed in §11.

## 6. There is No Representation of the dSA based on a One-Block Backbone

From CR's (2.2) it follows that $[J_z,[J_z,V_x]] = J_z^2V_x - 2J_zV_xJ_z + V_xJ_z^2 = V_x$. Using the explicit expression, (4.4e), for $J_z$ which applies for the HLA irrep $(A,B)$ this implies,

$$(V_x)_{ab,a'b'} = [(a+b)-(a'+b')]^2(V_x)_{ab,a'b'} \tag{6.1}$$

But the CRs (2.2) also require $[K_z,V_x]=0$ and using (4.4f) for $K_z$ this gives,

$$[(a-b)-(a'-b')](V_x)_{ab,a'b'} = 0 \tag{6.2}$$

The only way we can avoid concluding from (6.1) and (6.2) that $(V_x)_{ab,a'b'}=0$ is if both,

$$[(a+b)-(a'+b')] = \pm 1 \quad \text{and} \quad [(a-b)-(a'-b')] = 0 \tag{6.3}$$

But adding these implies $(a-a') = \pm 1/2$, which is impossible because the indices can differ only by an integer. Hence we conclude that $V_x = 0$. The other CRs in (2.2) then immediately imply that all the $V_\mu$ are zero, but the CRs (2.3) cannot then hold and so there is no such representation. QED.

## 7. The Minimum A and B in Any Representation is Zero

As a preliminary we derive expressions for the on-BD blocks of the commutator $[V_x,V_y]$. Note that due to the CRs (2.3), the on-BD block $[V_x,V_y]_{II}$ is required to be identically equal to $i(J_z)_{II}$, which is a diagonal matrix, given by (4.4e). Similarly, the on-BD block $[V_t,V_z]_{II}$ is required to be identically equal to $i(K_z)_{II}$, which is a diagonal matrix, given by (4.4f).

[In fact, all the commutators $[V_\alpha,V_\beta]$ must be block-diagonal, due to (2.3), and the requirement for their off-BD blocks to be identically zero is taken up in §8].

### 7.1 Expressions for the On-BD Components of $[V_x,V_y]$

The on-BD blocks of $[V_x,V_y]$ can be written,

$$[V_x,V_y]_{II} = V_{xIJ}V_{yJI} + V_{xIK}V_{yKI} + \cdots - (V_{yIJ}V_{xJI} + V_{yIK}V_{xKI} + \cdots) \tag{7.1}$$

where $J,K,\ldots$ are all the blocks connected to block $I$. We can write the contribution of each connected block $J,K,\ldots$ to (7.1) as,

$$t_{IJ}t_{JI}Z_{IJ} = i(V_{xIJ}V_{yJI} - V_{yIJ}V_{xJI}) \tag{7.2}$$

Here $t_{IJ},t_{JI}$ are the parameters introduced in Appendix A. We thus have,

$$i[V_x,V_y]_{II} = t_{IJ}t_{JI}Z_{IJ} + t_{IK}t_{KI}Z_{IK} + \cdots \tag{7.3a}$$

We have factored out the dependence on $t_{IJ},t_{JI}$ so that $Z_{IJ}$ has a unique expression independent of these unknown parameters. For the four possible cases of the connectivity between blocks $I$ and $J$ the expressions for $Z_{IJ}$ can be calculated using Table A.1 in Appendix A. The results are given in Table 1, below, in terms of $(A_I,B_I)$ and $(a_I,b_I)$. Care with signs is essential, the real quantity $i[V_x,V_y]_{II} \equiv -(J_z)_{II} \equiv -(a_I+b_I)$. In the latter expression and throughout this section we shall generally drop the factors of Kronecker deltas when referring to diagonal matrices, e.g., $\delta_{a_Ia_I'}\delta_{b_Ib_I'}$ ,which should strictly be present. These should be understood.

Similarly we can define $\tilde{Z}_{IJ}$ by,

$$[V_t, V_z]_{II} = t_{IJ}t_{JI}\tilde{Z}_{IJ} + t_{IK}t_{KI}\tilde{Z}_{IK} + \cdots \tag{7.3b}$$

Again using Table A.1 it is found that the $\tilde{Z}_{IJ}$ are equal to the $Z_{IJ}$ of Table 1 except with the sign of the terms in $a_I$ reversed. As $[V_t, V_z]_{II} = \mathrm{i}(K_z)_{II} = a_I - b_I$ is the same as $-(J_z)_{II}$ with the sign of $a_I$ reversed, satisfaction of this CR also follows. Satisfaction of the other CRs (2.3) is automatic as Table A.1 was derived on this basis. Note that the sign signature in Table 1 is identified with the signs in $A_I = A_J \pm \frac{1}{2}$ and $B_I = B_J \pm \frac{1}{2}$ in that order.

**Table 1: The $Z_{IJ}$ in (9.2) and (9.3a).**

| Connectivity $IJ$ | $Z_{IJ}$ |
|---|---|
| $++$ | $4[A_I b_I + B_I a_I]$ |
| $--$ | $-4[(A_I+1)b_I + (B_I+1)a_I]$ |
| $+-$ | $4[A_I b_I - (B_I+1)a_I]$ |
| $-+$ | $-4[(A_I+1)b_I - B_I a_I]$ |

### 7.2 The Minimum A and B in Any Representation is Zero

**Figure 1**

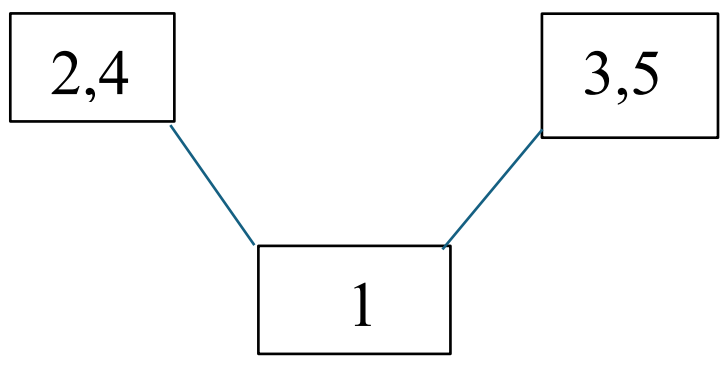

Proof: A block at minimum $A = A_{min}$, say $(A_{min}, B)$, can be connected only to blocks at $\left(A_{min} + \frac{1}{2}, B - \frac{1}{2}\right)$ and/or at $\left(A_{min} + \frac{1}{2}, B + \frac{1}{2}\right)$. There may be multiple copies of such blocks, say blocks 2 and 4 at the former and blocks 3 and 5 at the latter. This connectivity is illustrated in Figure 1 which plots $A$ vertically and $B$ horizontally. This type of connectivity diagram will be used extensively from this point. The lines between the boxes in Figure 1 denote connection, i.e., that the corresponding off-BD blocks of $V_\mu$ are non-zero.

The connections 12 and 14 are of type $-+$ whereas connections 13 and 15 are $--$. Using (7.3a) and Table 1 gives, (7.4)

$$\frac{i}{4}\left[V_x, V_y\right]_{11} = -(t_{12}t_{21} + t_{14}t_{41})[(A+1)b_1 - Ba_1] - (t_{13}t_{31} + t_{15}t_{51})[(A+1)b_1 + (B+1)a_1]$$

But from (A.21a,b) we have $t_{12}^* = t_{21}$, $t_{14}^* = t_{41}$, $t_{13}^* = -t_{31}$, $t_{15}^* = -t_{51}$, so (7.4) becomes,

$$\frac{i}{4}\left[V_x, V_y\right]_{11} = -(|t_{12}|^2 + |t_{14}|^2)[(A+1)b_1 - Ba_1] + (|t_{13}|^2 + |t_{15}|^2)[(A+1)b_1 + (B+1)a_1]$$

$$= b_1(A+1)[(|t_{13}|^2 + |t_{15}|^2) - (|t_{12}|^2 + |t_{14}|^2)] + a_1[B(|t_{12}|^2 + |t_{14}|^2) + (B+1)(|t_{13}|^2 + |t_{15}|^2)] \tag{7.5}$$

But CRs (2.3) require $\frac{i}{4}\left[V_x, V_y\right]_{11} \equiv -(J_z)_{11} = -(a_1 + b_1)$ and we see that this cannot be identical with (7.5), in which the coefficient of $a_1$ is non-negative, unless $A = A_{min} = 0$ so that $a_1 = 0$.

By symmetry the same must hold for the smallest $B$, i.e., $B_{min} = 0$. QED.

A further observation which follows from (7.5) is that any block $(A, B)$ which is connected only to a block (or blocks) "above" it, i.e., at $\left(A+\frac{1}{2}, B \pm \frac{1}{2}\right)$, must have $A=0$ if (7.5) is not to result in a contradiction.

By symmetry we see that any block $(A, B)$ which is connected only to a block (or blocks) "to the right" of it, i.e., at $\left(A \pm \frac{1}{2}, B+\frac{1}{2}\right)$, must have $B=0$ if it is not to result in a contradiction.

Putting those two results together means that any block connected only to blocks above and to its right must be the origin, $(0,0)$. This result also follows from (7.5) when blocks 2 and 4 are omitted.

## 8. Exclusion of Unique Non-Monotonic 3-Sequences

Considering a given sequence of three blocks, $IKJ$, whether monotonic or not, the corresponding off-BD block $IJ$ of $[V_x, V_y]$ will be denoted,

$$i[V_x, V_y]_{IJ} = t_{IK} t_{KJ} Z_{IKJ} \tag{8.1}$$

Here we have factored out the dependence on $t_{IK}$ and $t_{KJ}$ so that $Z_{IKJ}$ has a unique expression independent of these unknown parameters. For a given block $I$ there are 16 possible sequences of three blocks, $IKJ$, of which 4 are monotonic and 12 non-monotonic. Using Table A.1 the resulting $Z_{IKJ}$ can be evaluated in terms of $(A_I, B_I)$ and $(a_I, b_I)$. The results are given in Table 2. Note that due to the CRs (2.3), and the fact that $J_i$ are block-diagonal, the off-BD blocks $[V_x, V_y]_{IJ}$ $(I \neq J)$ are required to be identically zero.

Table 2 shows why monotonic sequences are special in this respect as they all have $Z_{IKJ}$ which are identically zero.

In contrast, if in any proposed representation there is a unique non-monotonic sequence of three blocks between any pair of blocks $I$ and $J$, the proposed representation cannot be valid because the corresponding $Z_{IKJ}$ is not identically zero and so $[V_x, V_y]_{IJ}$ fails to be identically zero. (The cases $A_I = 0$ and/or $B_I = 0$ may seem to be an exception to this for some sequences in Table 2, but zero values for $A_I$ and/or $B_I$ are incompatible with the sequence in those cases, so there is no exception. The same applies to some cases which would be identically zero if $A_I$ or $B_I$ were 1/2).

Note, however, that the stipulation that the non-monotonic 3-sequence be unique is essential to this conclusion. If there were two or more non-monotonic 3-sequences connecting blocks $I$ and $J$, say $IKJ$, $ILJ$, etc., then (8.1) would become instead,

$$i[V_x, V_y]_{IJ} = t_{IK} t_{KJ} Z_{IKJ} + t_{IL} t_{LJ} Z_{ILJ} + \cdots \tag{8.2}$$

This does not necessarily create a conflict with the requirement that $[V_x, V_y]_{IJ}$ be identically zero. This is because some of the expressions of Table 2 become linearly dependent upon others for specific values of $A_I$ and $B_I$, so that the RHS of (8.2) could be identically zero if the values of $t_{IK} t_{KJ}$ and $t_{IL} t_{LJ}$ were chosen appropriately. For example, if the blocks $I, K, J, L$ have $(A, B)$ values of $\left(1, \frac{1}{2}\right), \left(\frac{1}{2}, 0\right), \left(0, \frac{1}{2}\right), \left(\frac{1}{2}, 1\right)$ respectively, and if all compatible blocks are connected, then $IKJ$ has signature $++/+-$ and $ILJ$ has signature $+-/++$. Using Table 2, (8.2) then gives,

$$i[V_x, V_y]_{IJ} = (-2t_{IK}t_{KJ} + 6t_{IL}t_{LJ})\sqrt{1 - a_I^2} \tag{8.3}$$

With suitable choice of $t_{IK}, t_{KJ}, t_{IL}, t_{LJ}$ (8.3) can be identically zero. We will see later that this example does give rise to a valid irrep. This is why it is explicitly ***unique*** non-monotonic sequences which will be the focus of our attention, as these unambiguously indicate that a representation does not exist with the proposed connectivity.

We can also define $\tilde{Z}_{IKJ}$ via,

$$[V_t, V_z]_{IJ} = t_{IK}t_{KJ}\tilde{Z}_{IKJ} + t_{IL}t_{LJ}\tilde{Z}_{ILJ} + \cdots \tag{8.4}$$

These $\tilde{Z}_{IKJ}$ may be shown to equal the $Z_{IKJ}$ of Table 2 but with the signs of all terms involving $a_I$ reversed.

**Table 2: The $Z_{IKJ}$ in (8.1) and (8.2).** The stated sequence relates to $IK, KJ$.

| Sequence $IK, KJ$ | $Z_{IKJ}$ |
|---|---|
| $++,++$ | $0$ |
| $++,--$ | $4[A_I b_I + B_I a_I]$ |
| $++,+-$ | $-4B_I\sqrt{A_I^2 - a_I^2}$ |
| $++,-+$ | $-4A_I\sqrt{B_I^2 - b_I^2}$ |
| $--,++$ | $-4[(A_I+1)b_I + (B_I+1)a_I]$ |
| $--,--$ | $0$ |
| $--,+-$ | $-4(A_I+1)\sqrt{(B_I+1)^2 - b_I^2}$ |
| $--,-+$ | $-4(B_I+1)\sqrt{(A_I+1)^2 - a_I^2}$ |
| $+-,++$ | $4(B_I+1)\sqrt{A_I^2 - a_I^2}$ |
| $+-,--$ | $4A_I\sqrt{(B_I+1)^2 - b_I^2}$ |
| $+-,+-$ | $0$ |
| $+-,-+$ | $4[A_I b_I - (B_I+1)a_I]$ |
| $-+,++$ | $4(A_I+1)\sqrt{B_I^2 - b_I^2}$ |
| $-+,--$ | $4B_I\sqrt{(A_I+1)^2 - a_I^2}$ |
| $-+,+-$ | $4[B_I a_I - (A_I+1)b_I]$ |
| $-+,-+$ | $0$ |

## 9. Diagrammatic Classification of All Finite Irreps

In this section we present the classification scheme which we claim covers all the finite irreps (with Hermitian or anti-Hermitian generators). This done by exhaustive account of all backbones-plus-connectivity which respect the rules derived in §6, §7, §8 and Appendix A, namely,

- The minimum number of blocks in the backbone is 2;
- The smallest $A$ and the smallest $B$ in the backbone are both zero;
- There are no unique non-monotonic 3-sequences.
- Connection occurs only for $\left|A_P - A_Q\right| = \frac{1}{2}$ and $\left|B_P - B_Q\right| = \frac{1}{2}$.

To these we add the condition that no block occurs more than once in the backbone. The legitimacy of this assumption is justified in Appendix C.

The existence of the irreps thus classified is proved by explicit construction for two infinite classes in Appendix B. These classes are defined as the cases with backbones and connections which are fully monotonic. For the remaining cases, the existence of the irreps is proved in §10 for a few of the smallest dimensional instances by explicit construction. A proof that all such irreps exist is outstanding but seems likely.

Altogether there is an infinite discrete set of types of irrep, each one of which is an infinite discrete class of irreps. The irreps are thus uniquely categorised by two integers, written $m/n$, the meaning of which relates to the diagrammatic scheme, below. The relationship between this $m/n$ classification and the usual two parameter classification of the Casimir invariants (usually denoted $p$ and $q$) will be obtained in §11.

The diagrammatic scheme is based on plotting the position of the backbone blocks on an $A$ versus $B$ plot, together with connecting lines between these points (i.e., blocks) whose presence denotes connection is the sense defined in §5 (i.e., that the associated off-BD terms in $V_\mu$ are non-zero).

Note that the connection rule means that two blocks can be connected only if they are adjacent along a diagonal line of slope +1 or -1.

The only two-block irreps obeying the above rules are,

***Figure 2: The Two 2-Block Irreps (denoted 2/1 and 1/2 respectively).***

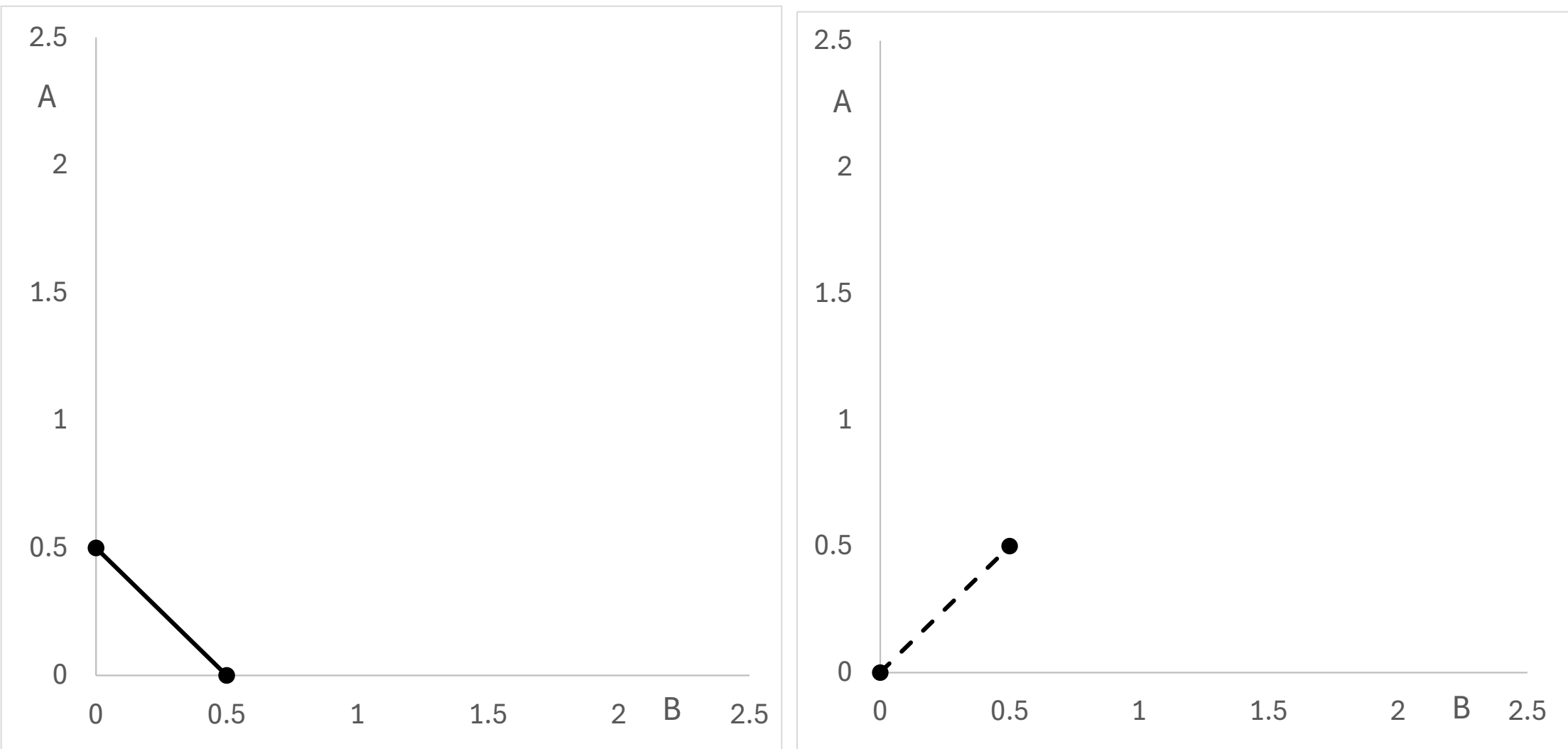

A non-monotonic 3-sequence is a "corner" in such a diagram. For example, both diagrams in Figure 3 contain a unique non-monotonic 3-sequence and therefore cannot be representations, despite obeying the other rules.

***Figure 3: Examples that are Not Representations***

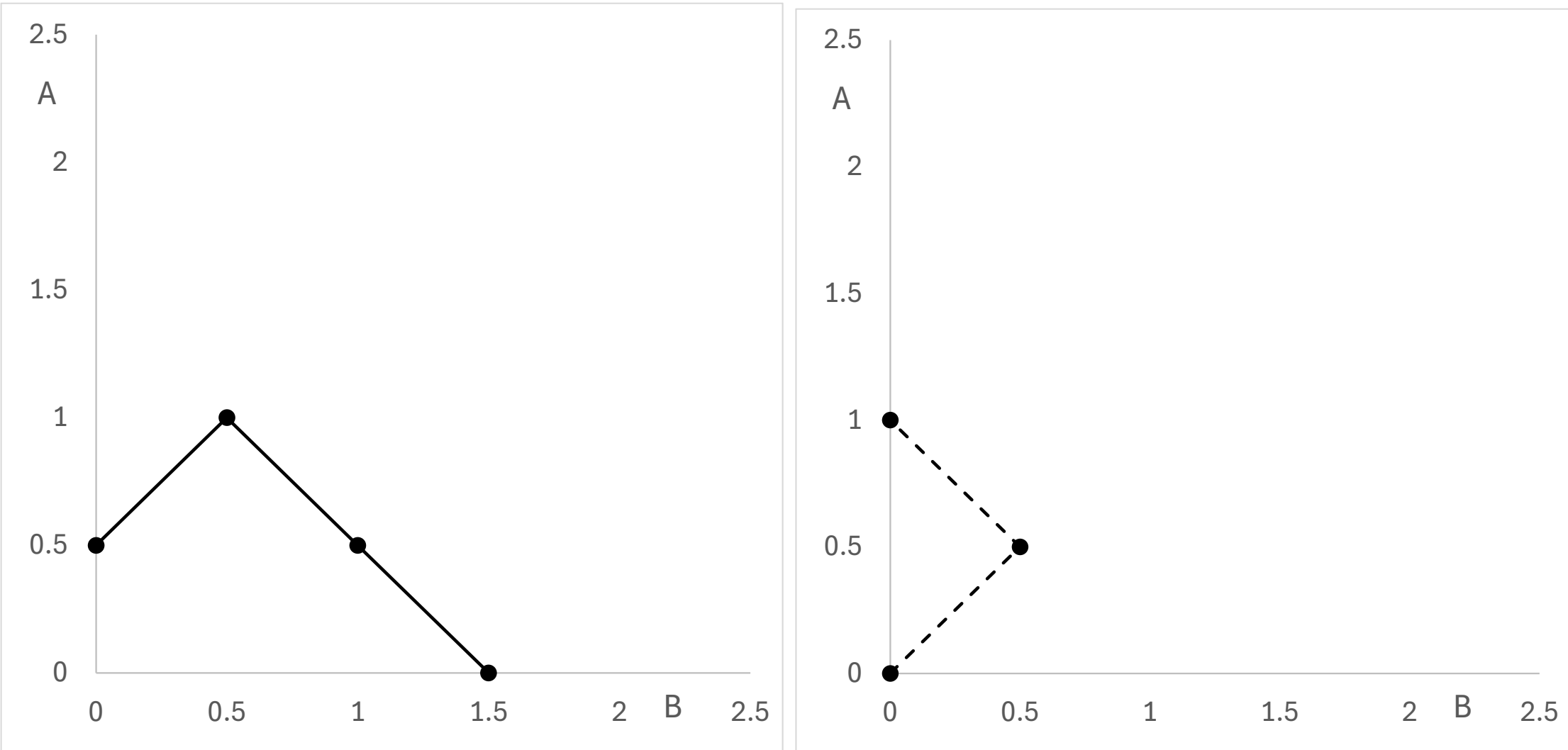

Note that in all these diagrams bold lines denote fermionic cases (half-integral spin) and dashed lines bosonic cases (integral spin).

The obvious way in which to avoid non-monotonic 3-sequences is to have diagrams in which the connections lie on a single straight line. The connectivity rules mean that these must be lines of gradient +1 or -1. Respecting also the rule that the minima of $A$ and $B$ are both zero gives us the first two infinite classes or irrep, illustrated in Figure 4.

***Figure 4: Illustrating the Two Infinite Classes of Fully Monotonic Irreps, denoted $\mathbf{m/1}$ and $\mathbf{1/n}$ respectively ($\mathbf{m, n} \in \{\mathbf{2, 3, 4, \ldots}\}$). Note that for $\mathbf{1/n}$ the block may terminate at any point for $\mathbf{n > 1}$.***

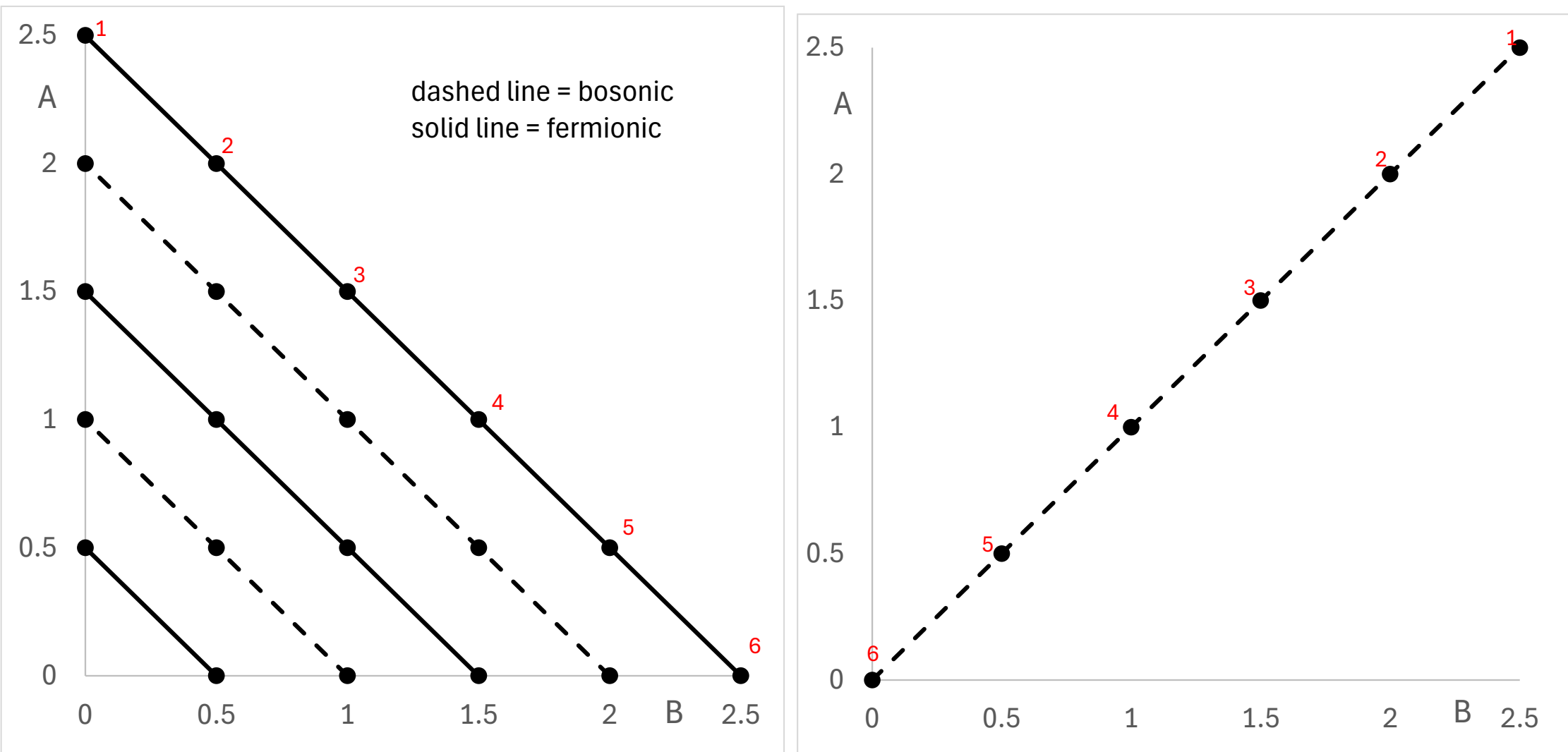

In the $m/n$ designation of the irreps, $m$ is the number of blocks in the -1 direction, whereas $n$ is the number of blocks in the +1 direction.

The existence of all the irreps in these two infinite classes is proved by explicit construction in Appendix B. This Appendix derives expressions for the parameters $t_{IJ}, t_{JI}$ which, together with Table A.1, fully specifies the representations of the $V_\mu$ (and the $J_i$ and $K_i$ are given by equs.(4.4a-f)). The two classes of irrep illustrated by Figure 4 have backbones as follows,

$$(A,A) \oplus \left(A-\tfrac{1}{2}, A-\tfrac{1}{2}\right) \oplus (A-1, A-1) \oplus \dots \oplus (0,0) \qquad (9.1a)$$

Or

$$(A,0) \oplus \left(A-\tfrac{1}{2}, \tfrac{1}{2}\right) \oplus (A-1,1) \oplus \dots \oplus (0,A) \qquad (9.1b)$$

in which there are at least two blocks.

Necessarily, $A = \frac{N-1}{2}$ where $N$ is the number of blocks in the backbone ($N \geq 2$). Hence, in order of ascending dimension, the first ten fully monotonic irreps' backbones are given in Table 3. The possible dimensions are $\frac{1}{6}N(N+1)(N+2)$ for Type $N/1$, whereas for type $1/N$ the dimension is $\frac{1}{6}N(N+1)(2N+1)$, where $N \in \{2,3,4,\dots\}$ is the number of blocks in the backbone in both cases.

**Table 3: The Backbones of the First Ten Fully Monotonic Irreps.** The second figure is the dimension of the representation.

| Rep | Dim | Backbone |
|---|---|---|
| 2/1 | 4 | $\left(\frac{1}{2},0\right) \oplus \left(0,\frac{1}{2}\right)$ |
| 1/2 | 5 | $\left(\frac{1}{2},\frac{1}{2}\right) \oplus (0,0)$ |
| 3/1 | 10 | $(1,0) \oplus \left(\frac{1}{2},\frac{1}{2}\right) \oplus (0,1)$ |
| 1/3 | 14 | $(1,1) \oplus \left(\frac{1}{2},\frac{1}{2}\right) \oplus (0,0)$ |
| 4/1 | 20 | $\left(\frac{3}{2},0\right) \oplus \left(1,\frac{1}{2}\right) \oplus \left(\frac{1}{2},1\right) \oplus \left(0,\frac{3}{2}\right)$ |
| 1/4 | 30 | $\left(\frac{3}{2},\frac{3}{2}\right) \oplus (1,1) \oplus \left(\frac{1}{2},\frac{1}{2}\right) \oplus (0,0)$ |
| 5/1 | 35 | $(2,0) \oplus \left(\frac{3}{2},\frac{1}{2}\right) \oplus (1,1) \oplus \left(\frac{1}{2},\frac{3}{2}\right) \oplus (0,2)$ |
| 1/5 | 55 | $(2,2) \oplus \left(\frac{3}{2},\frac{3}{2}\right) \oplus (1,1) \oplus \left(\frac{1}{2},\frac{1}{2}\right) \oplus (0,0)$ |
| 6/1 | 56 | $\left(\frac{5}{2},0\right) \oplus \left(2,\frac{1}{2}\right) \oplus \left(\frac{3}{2},1\right) \oplus \left(1,\frac{3}{2}\right) \oplus \left(\frac{1}{2},2\right) \oplus \left(0,\frac{5}{2}\right)$ |
| 1/6 | 91 | $\left(\frac{5}{2},\frac{5}{2}\right) \oplus (2,2) \oplus \left(\frac{3}{2},\frac{3}{2}\right) \oplus (1,1) \oplus \left(\frac{1}{2},\frac{1}{2}\right) \oplus (0,0)$ |

What remains are irreps which do contain non-monotonic 3-sequences (i.e., "corners") but that these are all non-unique. Assuming no duplicate blocks, the only means by which a "corner" comprising three non-monotonic blocks $IKJ$ may be made non-unique is to complete the square by including the block L with connectivity $ILJ$. The smallest dimension irrep including non-monotonic sequences is thus as illustrated in Figure 5 and denoted 2/2. This gives rise (we claim) to an infinite sequence of irreps, denoted $2/N$, for $N \in \{2,3,4,\dots\}$ of which the next in sequence, 2/3, is also shown in Figure 5.

***Figure 5: The Irreps 2/2 and 2/3, Illustrating the Series*** $\boldsymbol{2/N}$ ***for*** $\boldsymbol{N \in \{2,3,4,\dots\}}$

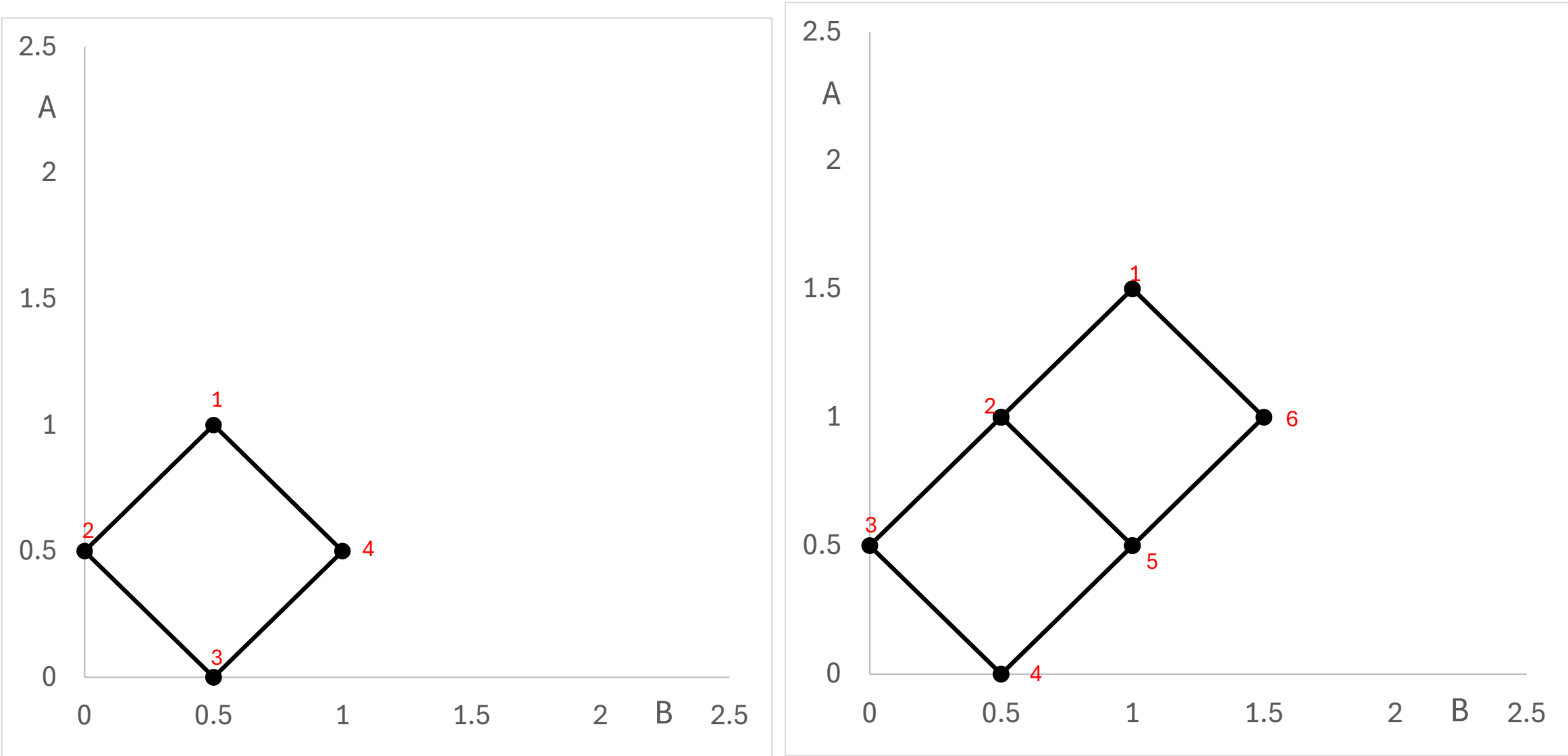

Similarly, Figure 6 shows the irreps 3/2 and 3/3, illustrating the next infinite class, $3/N$ for $N \in \{2,3,4,\dots\}$.

***Figure 6: The Irreps 3/2 and 3/3, Illustrating the Series*** $\boldsymbol{3/N}$ ***for*** $\boldsymbol{N \in \{2,3,4,\dots\}}$

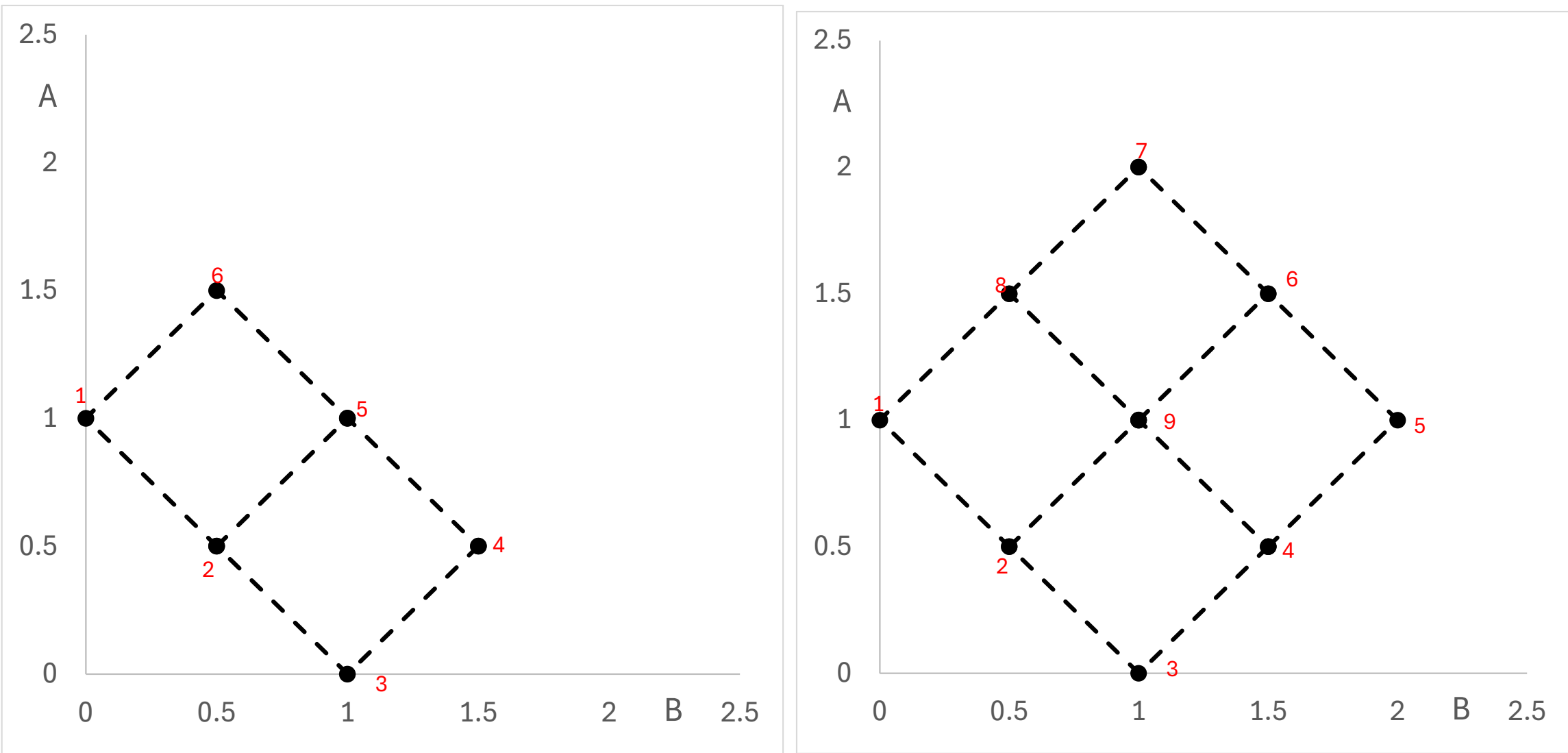

Similarly, Figure 7 shows the irreps 4/2 and 4/3, illustrating the next infinite class, $4/N$ for $N \in \{2,3,4,\dots\}$, and so on for an infinite number of infinite classes.

*Figure 7: The Irreps 4/2 and 4/3, Illustrating the Series $4/N$ for $N \in \{2, 3, 4, \dots\}$*

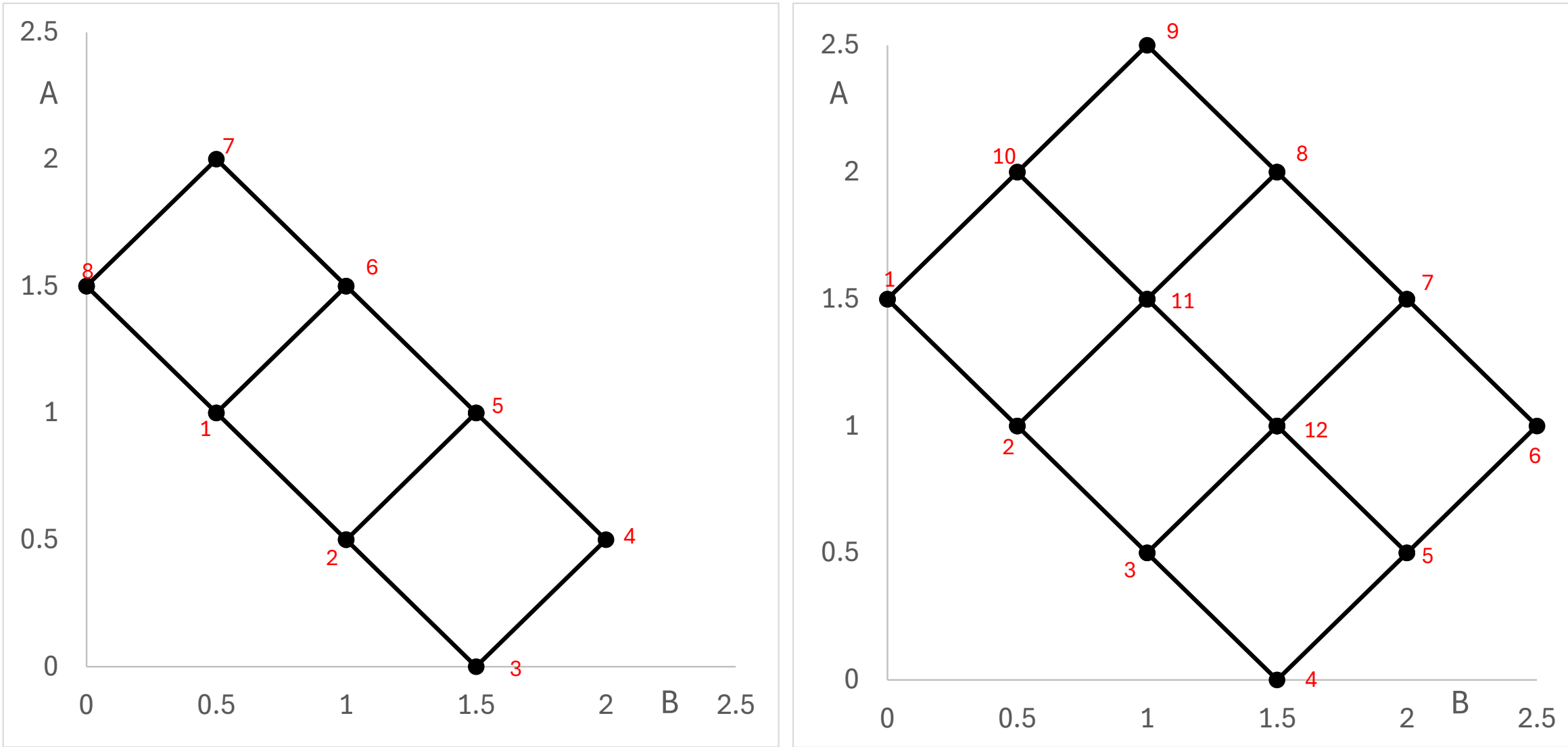

Note that the omission of any single connection in Figures 5, 6, 7 would inevitably leave at least one pair of blocks with a unique non-monotonic 3-sequence, and hence would not be a representation. The irreps described above are therefore exhaustive; no other finite irreps of the dSA with Hermitian or anti-Hermitian generators exists.

The existence of the irreps illustrated in Figures 5, 6 and 7 have been proved by explicit construction (see §10). We claim that the irreps $m/n$ for all $m, n \in \{1,2,3,4, \dots\}$, excluding only $m = n = 1$, all exist, though the proof is outstanding.

## 10. Construction of the Irreps Corresponding to Figures 5, 6 and 7

The full set of representation matrices are defined by equs.(4.4a-f) for $J_i$ and $K_i$, whilst the $V_\mu$ are defined by the combination of, (i) the backbone and its connectivity, which are given by Figures 5-7, (ii) the generic form of the $U_{IJ}$ matrices (see equ.(A.19) and Table A.1), and, (iii) specifying the values of the $t_{IJ}$ parameters which occur in equ.(A.19). The latter are given for the irreps of Figures 5-7 in Table 4 where the block numbers are defined in those Figures. The dimensions of the irreps are given in each case.

***Tables 4: The $t_{IJ}$ Parameters Defining the Irrep Matrices for Figs 5-7. In all cases $t_{IJ} = t_{JI}$***

***Table 4a: Irrep 2/2 (dimension 16)***

| | | | |
|---|---|---|---|
| $t_{12}$ | $\frac{\sqrt{5}}{4}i$ | $t_{34}$ | $\frac{\sqrt{5}}{4}i$ |
| $t_{23}$ | $\frac{3}{4}$ | $t_{14}$ | $\frac{1}{4}$ |

***Table 4b: Irrep 2/3 (dimension 40)***

| | | | |
|---|---|---|---|
| $t_{12}$ | $\sqrt{\frac{7}{36}}\ i$ | $t_{56}$ | $\sqrt{\frac{7}{36}}\ i$ |
| $t_{23}$ | $\frac{\sqrt{3}}{2}i$ | $t_{16}$ | $\frac{1}{6}$ |
| $t_{34}$ | 1 | $t_{25}$ | $\frac{1}{3}$ |
| $t_{45}$ | $\frac{\sqrt{3}}{2}i$ | | |

***Table 4c: Irrep 2/4 (dimension 80)***

| | | | | | |
|---|---|---|---|---|---|
| $t_{12}$ | $\frac{2}{3}i$ | $t_{56}$ | $\frac{2}{3}i$ | $t_{81}$ | $\frac{3}{8}i$ |
| $t_{23}$ | $\sqrt{\frac{63}{48}}\ i$ | $t_{67}$ | $\sqrt{\frac{27}{192}}\ i$ | $t_{25}$ | $\frac{5}{12}$ |
| $t_{34}$ | $\sqrt{\frac{75}{48}}$ | $t_{78}$ | $\frac{1}{8}$ | $t_{16}$ | $\frac{5}{24}$ |
| $t_{45}$ | $\sqrt{\frac{63}{48}}\ i$ | | | | |

***Table 4d: Irrep 3/2 (dimension 35)***

| | | | |
|---|---|---|---|
| $t_{12}$ | $\frac{1}{\sqrt{2}}$ | $t_{56}$ | $\frac{1}{\sqrt{12}}$ |
| $t_{23}$ | $\frac{1}{\sqrt{2}}$ | $t_{16}$ | $\frac{i}{2}$ |
| $t_{34}$ | $\frac{i}{2}$ | $t_{25}$ | $\frac{i}{\sqrt{6}}$ |
| $t_{45}$ | $\frac{1}{\sqrt{12}}$ | | |

***Table 4e: Irrep 3/3 (dimension 81)***

| | | | | | |
|---|---|---|---|---|---|
| $t_{12}$ | $\sqrt{\frac{5}{6}}$ | $t_{56}$ | $\frac{1}{\sqrt{24}}$ | $t_{29}$ | $\sqrt{\frac{7}{18}}\ i$ |
| $t_{23}$ | $\sqrt{\frac{5}{6}}$ | $t_{67}$ | $\frac{1}{\sqrt{24}}$ | $t_{49}$ | $\frac{\sqrt{5}}{6}$ |
| $t_{34}$ | $\sqrt{\frac{7}{12}}\ i$ | $t_{78}$ | $\sqrt{\frac{1}{6}}\ i$ | $t_{69}$ | $\frac{\sqrt{5}}{6}\ i$ |
| $t_{45}$ | $\sqrt{\frac{1}{6}}\ i$ | $t_{81}$ | $\sqrt{\frac{7}{12}}\ i$ | $t_{89}$ | $\frac{\sqrt{5}}{6}$ |

***Table 4f: Irrep 4/2 (dimension 64)***

| | | | | | |
|---|---|---|---|---|---|
| $t_{12}$ | $\frac{2}{3}$ | $t_{56}$ | $\frac{1}{3}$ | $t_{16}$ | $\sqrt{\frac{35}{288}}\ i$ |
| $t_{23}$ | $\sqrt{\frac{15}{32}}$ | $t_{67}$ | $\sqrt{\frac{3}{32}}$ | $t_{25}$ | $\sqrt{\frac{35}{288}}\ i$ |
| $t_{34}$ | $\sqrt{\frac{7}{32}}\ i$ | $t_{78}$ | $\sqrt{\frac{7}{32}}\ i$ | | |
| $t_{45}$ | $\sqrt{\frac{3}{32}}$ | $t_{81}$ | $\sqrt{\frac{15}{32}}$ | | |

*Table 4g: Irrep 4/3 (dimension 140)*

| | | | | | | | |
|---|---|---|---|---|---|---|---|
| $t_{12}$ | $\frac{\sqrt{3}}{2}$ | $t_{67}$ | $\sqrt{\frac{1}{20}}$ | $t_{10,11}$ | $\sqrt{\frac{3}{20}}$ | $t_{3,12}$ | $\sqrt{\frac{5}{18}}\,i$ |
| $t_{23}$ | $\frac{5}{6}$ | $t_{78}$ | $\sqrt{\frac{5}{80}}$ | $t_{8,11}$ | $\sqrt{\frac{9}{80}}\,i$ | $t_{5,12}$ | $\sqrt{\frac{3}{20}}$ |
| $t_{34}$ | $\frac{\sqrt{3}}{2}$ | $t_{89}$ | $\sqrt{\frac{1}{20}}$ | $t_{2,11}$ | $\sqrt{\frac{5}{18}}\,i$ | | |
| $t_{45}$ | $\frac{i}{\sqrt{2}}$ | $t_{9,10}$ | $\sqrt{\frac{3}{20}}\,i$ | $t_{11,12}$ | $\frac{5}{12}$ | | |
| $t_{56}$ | $\sqrt{\frac{3}{20}}\,i$ | $t_{10,1}$ | $\frac{i}{\sqrt{2}}$ | $t_{7,12}$ | $\sqrt{\frac{9}{80}}\,i$ | | |

There are 17 irreps with dimension less than 105. These are listed in Table 5.

***Table 5: All the Irreps with Dimension*** $< \mathbf{105}$

| irrep | 2/1 | 1/2 | 3/1 | 1/3 | 2/2 | 4/1 | 1/4 | 5/1 | 3/2 | 2/3 | 1/5 | 6/1 | 4/2 | 2/4 | 3/3 | 7/1 | 1/6 |
|---|---|---|---|---|---|---|---|---|---|---|---|---|---|---|---|---|---|
| $N_{dim}$ | 4 | 5 | 10 | 14 | 16 | 20 | 30 | 35 | 35 | 40 | 55 | 56 | 64 | 80 | 81 | 84 | 91 |

## 11. The Casimir Matrices

The Casimir operators are,

$$C1 = V_t^2 + K^2 - J^2 - V^2 \tag{11.1}$$

$$C2 = (\bar{K}\cdot\bar{J})^2 - (\bar{V}\cdot\bar{J})^2 + \mathcal{J}^2 \tag{11.2}$$

where $V^2 = V_x^2 + V_y^2 + V_z^2$, etc., and,

$$\bar{\mathcal{J}} = V_t\bar{J} + \bar{K}\times\bar{V} \tag{11.3}$$

The Casimir operators commute with all the generators of the Lie algebra and hence their representations within an irrep equals a unit matrix times a constant, the Casimir invariants. The possible values of the Casimir invariants for the irreps can be derived by explicit substitution of the complete representation matrices derived above. For C1 this is facilitated by the following algebraically equivalent expression, (11.4)

$$C1 = K_z^2 - J_z^2 + \frac{1}{2}\{(K_+K_- + K_-K_+) - (J_+J_- + J_-J_+)\} - 2\{(V_+V_- + V_-V_+) + (W_+W_- + W_-W_+)\}$$

Each of the terms in (11.4) can now be evaluated using (4.4a-f) and Tables A.1, B.1 and 4. Doing so confirms that the Casimir matrices are indeed multiples of the unit matrix and hence confirms that the representations as described above are irreps (and the only irreps). Direct sums of multiple copies of the same irrep would also produce Casimir matrices which were multiples of the unit matrix, but their reducible nature is then betrayed by the presence of duplicates.

As an example, consider the fully monotonic irreps illustrated by Figure 4 and Table 3. To find the Casimir invariants it suffices to evaluate block 11. A diagonal matrix results, as it must, and we find the Casimir invariants to be

For Type 1/n: $$-C1 = 4A_1(A_1+1) + 8t_{12}^2 A_1^2 \quad (11.5a)$$

For Type m/1: $$-C1 = 2A_1(A_1+1) + 8t_{12}^2 A_1 \quad (11.5b)$$

where $A_1$ is the greatest $A$. As noted in §1, work on the infinite dimensional unitary representations has resulted in $C1$ being expressed in terms of two parameters, which, in the case of the discrete representations, take positive integral or half-integral values. Various conventions have been used. In that of Dixmier and Enayati et al, the expression is,

$$-C1 = p(p+1) + (q+1)(q-2) \quad (11.6)$$

in which $p$ takes values from the infinite range $p \in \left\{\frac{1}{2}, 1, \frac{3}{2}, 2, \frac{5}{2} \dots\right\}$ and $q$ takes values from the finite range $q \in \left\{\frac{1}{2}, 1, \frac{3}{2}, 2, \frac{5}{2} \dots p\right\}$. (Infinite dimensional unitary representations are also found with the same expression for $C1$ but in which $q$ takes values from a continuum range).

Recall that for irreps Type n/1 we have $t_{21} = t_{12} = \frac{1}{2}$, whereas for irreps Type 1/n Table B.1 gives us $t_{12} = -t_{21} = \frac{1}{2\sqrt{A_1}}$. We thus find,

For Type 1/n: $$-C1 = 2A_1(2A_1+3) = p(p+1) - 2 \quad (11.7a)$$

where, $$p = 2A_1 + 1 \text{ and } q = 1 \text{ and } p \in \{2, 3, 4 \dots\} \quad (11.7b)$$

Here the value of $q$ assigned makes (13.7a) consistent with the infinite dimensional case, (13.6). The infinite dimensional unitary representations also provide a parametrised expression for $C2$,

$$-C2 = p(p+1)q(q-1) \quad (11.8)$$

This expression implies that $C2 = 0$ for Type 1/n for which $q = 1$. This is confirmed by direct computation.

Expressions (11.6) and (11.8) for the two Casimir invariants are also found to apply for Type m/1 but now with $q = p$, thus,

For Type m/1: $$-C1 = 2A_1(A_1+2) = 2(p^2-1) \quad (11.9a)$$

where, $$p = q = A_1 + 1 \text{ and } p \in \left\{\frac{3}{2}, 2, \frac{5}{2}, 3 \dots\right\} \quad (11.9b)$$

For the irreps of Type $m/n$ with $m, n > 1$ direct evaluation shows the same expressions, (11.6) and (11.8), apply. Overall, the Casimir invariants for the complete set of finite irreps, with Hermitian or anti-Hermitian generators are given by $p \in \left\{\frac{3}{2}, 2, \frac{5}{2}, 3 \dots\right\}$ and, for half-integral $p$ then $q$ is also half-integral with $q \in \left\{\frac{3}{2}, \frac{5}{2}, \dots p\right\}$ whereas for integral $p$ then $q$ is also integral with $q \in \{1,2, \dots p\}$.

We find that the $p, q$ parameters can be expressed simply in terms of the $m/n$ which characterise the irrep. These relations are,

$$q = \frac{m+1}{2} \quad \text{and} \quad p = n + \frac{m-1}{2} = n - 1 + q \tag{11.10}$$

Hence, the complete set of irreps is specified simply by all $m/n$ with $m, n \in \{1,2,3, \dots\}$ excepting only $m = n = 1$ and then the Casimir invariants are found immediately from,

$$-C1 = (n-1)(m+n+1) + \frac{1}{2}(m+1)^2 - 2 \tag{11.11a}$$

$$-C2 = \left(n^2 + \frac{m^2-1}{4} + mn\right)\left(\frac{m^2-1}{4}\right) \tag{11.11b}$$

Table 6 summarises the values taken by the Casimir invariants for all 21 irreps with dimension not exceeding 140.

**Table 6: The Values of the Casimir Scalars for the 21 Irreps up to Dimension 140**.

| Irrep | $N_{dim}$ | p | q | -Casimir1 | -Casimir2 |
|---|---|---|---|---|---|
| 2/1 | 4 | 1.5 | 1.5 | 5/2 | 45 /16 |
| 3/1 | 10 | 2 | 2 | 6 | 12 |
| 4/1 | 20 | 2.5 | 2.5 | 21/2 | 525 /16 |
| 5/1 | 35 | 3 | 3 | 16 | 72 |
| 6/1 | 56 | 3.5 | 3.5 | 45/2 | 2205/16 |
| 7/1 | 84 | 4 | 4 | 30 | 240 |
| 8/1 | 120 | 4.5 | 4.5 | 77/2 | 6237/16 |
| 1/2 | 5 | 2 | 1 | 4 | 0 |
| 1/3 | 14 | 3 | 1 | 10 | 0 |
| 1/4 | 30 | 4 | 1 | 18 | 0 |
| 1/5 | 55 | 5 | 1 | 28 | 0 |
| 1/6 | 91 | 6 | 1 | 40 | 0 |
| 1/7 | 140 | 7 | 1 | 54 | 0 |
| 2/2 | 16 | 2.5 | 1.5 | 15/2 | 105/16 |
| 2/3 | 40 | 3.5 | 1.5 | 29/2 | 189/16 |
| 2/4 | 80 | 4.5 | 1.5 | 47/2 | 297/16 |
| 3/2 | 35 | 3 | 2 | 12 | 24 |
| 3/3 | 81 | 4 | 2 | 20 | 40 |
| 4/2 | 64 | 3.5 | 2.5 | 35/2 | 945/16 |
| 4/3 | 140 | 4.5 | 2.5 | 53/2 | 1485/16 |
| 5/2 | 105 | 4 | 3 | 24 | 120 |

## 14. Conclusion

The complete set of all the finite irreps of the de Sitter Lie algebra with either Hermitian or anti-Hermitian generators have been characterised by a diagrammatic depiction conveniently labelled by two positive integers $m/n$ where $m, n \in \{1,2,3, \dots\}$ excepting only $m = n = 1$. The matrix elements have been specified, being defined by (4.4a-f) and Tables 4, A.1, B.1. The two Casimir invariants have been given explicitly in terms of the two positive integers, $m, n$, equs.(11.11a,b), and shown to reproduce the result familiar from the study of the infinite dimensional, unitary representations $-C1 = p(p+1) + (q+1)(q-2)$ and $-C2 = p(p+1)q(q-1)$ but with a reduced set of possible $p, q$ values.

## Appendix A: The Representations of $V_\mu$ and Compatibility between Blocks

This Appendix generalises to an arbitrary number of blocks the derivation by Shurtleff, Ref.[5], of the conditions under which block-components of $V_\mu$ can be non-zero. Note that Ref.[5] addressed the Poincaré group, not the dSG, but the derivations summarised in this section require only the CRs (2.1) and (2.2), which are common between the two algebras and hence Ref.[5] is applicable. Satisfaction of CRs (2.3) will be addressed in the main text and in Appendix B.

Considering for illustration the case of a backbone with four blocks, the representation of the $J_i$ and the $K_i$ can be written in block-matrix form as,

$$\begin{pmatrix} J_i^{(1)} & 0 & 0 & 0 \\ 0 & J_i^{(2)} & 0 & 0 \\ 0 & 0 & J_i^{(3)} & 0 \\ 0 & 0 & 0 & J_i^{(4)} \end{pmatrix} \quad \text{and} \quad \begin{pmatrix} K_i^{(1)} & 0 & 0 & 0 \\ 0 & K_i^{(2)} & 0 & 0 \\ 0 & 0 & K_i^{(3)} & 0 \\ 0 & 0 & 0 & K_i^{(4)} \end{pmatrix} \tag{A.1}$$

where $J_i^{(1)}$ is shorthand for $J_i^{(A_1,B_1)}$, $J_i^{(2)}$ is shorthand for $J_i^{(A_2,B_2)}$, etc., and similarly for the $K_i^{(P)}$. (We shall used capital Latin letters to denote block numbers). Note that $J_i^{(A_1,B_1)}$ is a matrix of dimension $(2A_1+1)(2B_1+1)$, etc. We initially write the rep of $V_\mu$ as a full matrix whose $PQ$ block is $V_{\mu PQ}$, each one of which is, in general, a rectangular matrix with $(2A_P+1)(2B_P+1)$ rows and $(2A_Q+1)(2B_Q+1)$ columns. The $PQ$ block of the CR $[J_i, V_\mu]$ is thus,

$$[J_i, V_\mu]_{PQ} = J_i^{(P)} V_{\mu PQ} - V_{\mu PQ} J_i^{(Q)} \tag{A.2}$$

with a similar expression for $[K_i, V_\mu]_{PQ}$. But on-block-diagonal (on-BD) blocks then do not mix blocks, i.e., $[J_i, V_\mu]_{PP} = J_i^{(P)} V_{\mu PP} - V_{\mu PP} J_i^{(Q)}$ involves only the on-BD block $PP$. But we have already shown in §6 that the only solution for the $V_{\mu PP}$ is that they are zero. Hence we conclude that the $V_\mu$ matrices have zero blocks on the BD.

We now present the proof that the off-BD block $V_{\mu PQ}$ between blocks $P$ and $Q$ can be non-zero only if $|A_P - A_Q| = \frac{1}{2}$ and $|B_P - B_Q| = \frac{1}{2}$.

Using the CRs (2.2) we have $[J_z, [J_z, V_x]] = J_z^2 V_x - 2J_z V_x J_z + V_x J_z^2 = V_x$. Considering the off-BD blocks 12 and 21, this decouples into two parts,

$$V_{x12} = \left(J_z^{A,B}\right)^2 V_{x12} - 2J_z^{A,B} V_{x12} J_z^{C,D} + V_{x12}\left(J_z^{C,D}\right)^2 \tag{A.3a}$$

$$V_{x21} = \left(J_z^{C,D}\right)^2 V_{x21} - 2J_z^{C,D} V_{x21} J_z^{A,B} + V_{x21}\left(J_z^{A,B}\right)^2 \tag{A.3b}$$

where we have taken blocks 1 and 2 to have HLA irreps $(A,B)$ and $(C,D)$ respectively. Using the explicit expression, (4.4e), for $J_z$ for both the $(A,B)$ and $(C,D)$ parts we find, from (A.3a),

$$(V_{x12})_{ab,cd} = [(a+b)-(c+d)]^2 (V_{x12})_{ab,cd} \tag{A.4}$$

where $c$ takes the integral or half-integral values $-C, -C+1 \dots . C$, and similarly for $d$ in terms of $D$. From (2.2), the CR $[K_z, V_x] = 0$ gives in similar fashion,

$$[(a-b)-(c-d)](V_{x12})_{ab,cd} = 0 \tag{A.5}$$

To avoid (A.4) and (A.5) implying $(V_{x12})_{ab,cd} = 0$ (and by immediate implication of the CRs, (2.2), that all the $V_\mu = 0$) we require,

$$[(a+b)-(c+d)] = \pm 1 \quad \text{and} \quad [(a-b)-(c-d)] = 0 \tag{A.6}$$

Adding the two equations implies $a-c = b-d = \pm 1/2$. This is possible only if $A$ is an integer and $C$ is a half-integer, or vice-versa, and similarly for $B$ and $D$. The conditions (A.6) establish the only components of $(V_{x12})_{ab,cd}$ which can be non-zero, namely those with either $a = c + \frac{1}{2}$ and $b = d + \frac{1}{2}$ or $a = c - \frac{1}{2}$ and $b = d - \frac{1}{2}$, and hence we can write,

$$(V_{x12})_{ab,cd} = t^{12}_{ab}\delta_{a,c+\frac{1}{2}}\delta_{b,d+\frac{1}{2}} + u^{12}_{ab}\delta_{a,c-\frac{1}{2}}\delta_{b,d-\frac{1}{2}} \tag{A.7}$$

By using the CR $[J_z, V_x] = iV_y$ we find that $i\left(V_{y12}\right)_{ab,cd} = [a+b-c-d](V_{x12})_{ab,cd}$ and hence,

$$i\left(V_{y12}\right)_{ab,cd} = t^{12}_{ab}\delta_{a,c+\frac{1}{2}}\delta_{b,d+\frac{1}{2}} - u^{12}_{ab}\delta_{a,c-\frac{1}{2}}\delta_{b,d-\frac{1}{2}} \tag{A.8}$$

Defining $V_{\pm 12} = \left(V_{x12} \pm iV_{y12}\right)/2$, (A.7,8) simplify to,

$$(V_{+12})_{ab,cd} = t^{12}_{ab}\delta_{a,c+\frac{1}{2}}\delta_{b,d+\frac{1}{2}} \qquad (V_{-12})_{ab,cd} = u^{12}_{ab}\delta_{a,c-\frac{1}{2}}\delta_{b,d-\frac{1}{2}} \tag{A.9a}$$

Symmetry means we can obtain the corresponding expressions for $V_{\pm 21}$ by the replacements $a \leftrightarrow c$, $b \leftrightarrow d$, giving,

$$(V_{+21})_{cd,ab} = t^{21}_{cd}\delta_{c,a+\frac{1}{2}}\delta_{d,b+\frac{1}{2}} \qquad (V_{-21})_{cd,ab} = u^{21}_{cd}\delta_{c,a-\frac{1}{2}}\delta_{d,b-\frac{1}{2}} \tag{A.9b}$$

So far we have used only the CRs involving $J_z$. Defining $J_\pm = J_x \pm iJ_y$ the CRs, (2.2), involving the $x$ and $y$ components can be simplified to,

$$[J_\pm, V_\pm] = 0 \qquad \text{and} \qquad [J_\pm, V_\mp] = \pm V_z \tag{A.10}$$

In terms of block components we can write,

$$([J_+, V_+]_{12})_{ab,cd} = \left(J^{A,B}_+\right)_{ab,a'b'}(V_{+12})_{a'b',cd} - (V_{+12})_{ab,c'd'}\left(J^{C,D}_+\right)_{c'd',cd} \tag{A.11}$$

where the dashed indices are summed. Using the explicit expressions (4.4a-f) together with (A.9a) this zero commutator gives,

$$\left(r^A_{a-1}t^{12}_{a-1,b} - r^C_c t^{12}_{ab}\right)\delta_{a,c+\frac{3}{2}}\delta_{b,d+\frac{1}{2}} = 0 \tag{A.12a}$$

and,
$$\left(r^B_{b-1}t^{12}_{a,b-1} - r^D_d t^{12}_{a,b}\right)\delta_{a,c+\frac{1}{2}}\delta_{b,d+\frac{3}{2}} = 0 \tag{A.12b}$$

In the same manner $([J_-, V_-]_{12})_{ab,cd} = 0$, $([J_+, V_+]_{21})_{cd,ab} = 0$ and $([J_-, V_-]_{21})_{cd,ab} = 0$ give respectively,

$$\left(s^A_{a+1}u^{12}_{a+1,b} - s^C_c u^{12}_{ab}\right)\delta_{a,c-\frac{3}{2}}\delta_{b,d-\frac{1}{2}} = 0 \tag{A.12c}$$

and,
$$\left(s^B_{b+1}u^{12}_{a,b+1} - s^D_d u^{12}_{a,b}\right)\delta_{a,c-\frac{1}{2}}\delta_{b,d-\frac{3}{2}} = 0 \tag{A.12d}$$

$$\left(r_{c-1}^{C}t_{c-1,d}^{21} - r_{a}^{A}t_{cd}^{21}\right)\delta_{c,a+\frac{3}{2}}\delta_{d,b+\frac{1}{2}} = 0 \qquad \text{(A.12e)}$$

and,

$$\left(r_{d-1}^{D}t_{c,d-1}^{21} - r_{b}^{B}t_{cd}^{21}\right)\delta_{c,a+\frac{1}{2}}\delta_{d,b+\frac{3}{2}} = 0 \qquad \text{(A.12f)}$$

$$\left(s_{c+1}^{C}u_{c+1,d}^{21} - s_{a}^{A}u_{cd}^{21}\right)\delta_{c,a-\frac{3}{2}}\delta_{d,b-\frac{1}{2}} = 0 \qquad \text{(A.12g)}$$

and,

$$\left(s_{d+1}^{D}u_{c,d+1}^{21} - s_{b}^{B}u_{cd}^{21}\right)\delta_{c,a-\frac{1}{2}}\delta_{d,b-\frac{3}{2}} = 0 \qquad \text{(A.12h)}$$

Consider (A.12a). So long as we choose $b$ and $d$ such that $b = d + \frac{1}{2}$ then (A.12a) implies, for any $a$ and $b$ compatible with this and with $a = c + \frac{3}{2}$,

$$r_{C}^{C}t_{c+\frac{3}{2},b}^{12} = r_{c+\frac{1}{2}}^{A}t_{c+\frac{1}{2},b}^{12} \qquad \text{(A.13)}$$

For any given $b$, (A.13) is a recursion formula for $t_{ab}^{12}$ as it gives each $t_{ab}^{12}$ in terms of the one before, i.e., $t_{a-1,b}^{12}$. Providing that $C \geq A + \frac{3}{2}$ we can choose $a = -A$ in (A.13) and still have a $c$ compatible with $a = c + \frac{3}{2}$. With $a = -A$, (A.13) gives,

$$r_{-A-\frac{3}{2}}^{C}t_{-A,b}^{12} = r_{-A-1}^{A}t_{-A-1,b}^{12} = 0 \qquad \text{(A.14)}$$

because (4.2d) gives $r_{-A-1}^{A} = 0$. If we assume this implies $t_{-A,b}^{12} = 0$ then with $a = -A + 1$, (A.13) gives,

$$r_{-A-\frac{1}{2}}^{C}t_{-A+1,b}^{12} = r_{-A}^{A}t_{-A,b}^{12} = 0 \qquad \text{(A.15)}$$

So $t_{-A+1,b}^{12} = 0$ also, and so on, finding that all the $t_{ab}^{12}$ are zero. The only way of avoiding this conclusion is if, in (A.14), $r_{-A-\frac{3}{2}}^{C} = \sqrt{\left(C + A + \frac{3}{2}\right)\left(C - A - \frac{1}{2}\right)} = 0$, i.e., $C = A + \frac{1}{2}$. But this contradicts the initial assumption that $C \geq A + \frac{3}{2}$. So we conclude that this initial assumption leads to all the $t_{ab}^{12}$ being zero.

The same assumption, $C \geq A + \frac{3}{2}$, can be used in the same way to deduce from (A.12e) that all the $t_{cd}^{21}$ are zero. And again, the assumption, $C \geq A + \frac{3}{2}$ can be used in the same way to deduce from (A.12c) and (A.12g) that all the $u_{ab}^{12}$ and all the $u_{ab}^{21}$ are zero. Hence, $C \geq A + \frac{3}{2}$ necessarily leads to $V_{\mu 12} = V_{\mu 21} = 0$.

The symmetry amongst the Equs.(A.12a-h) immediately implies that this conclusion will follow also for $A \geq C + \frac{3}{2}$. To derive this explicitly, note that we can now assume $c = C$ in (A.13) whilst being compatible with $a = c + \frac{3}{2}$. But this gives,

$$r_{C+\frac{1}{2}}^{A}t_{C+\frac{1}{2},b}^{12} = r_{C}^{C}t_{C+\frac{3}{2},b}^{12} = 0 \qquad \text{(A.16)}$$

because (4.2d) gives $r_{C}^{C} = 0$. And for $c = C - 1$ we get,

$$r_{C-\frac{1}{2}}^{A}t_{C-\frac{1}{2},b}^{12} = r_{C-1}^{C}t_{C+\frac{1}{2},b}^{12} = 0 \qquad \text{(A.17)}$$

and so on. If (C22) is taken to imply $t^{12}_{C+\frac{1}{2},b} = 0$ then this sequence gives all the $t^{12}_{ab} = 0$. To avoid this we require $r^{A}_{C+\frac{1}{2}} = 0 = \sqrt{\left(A - C - \frac{1}{2}\right)\left(A + C + \frac{3}{2}\right)}$ and hence $A = C + \frac{1}{2}$ which contradicts the initial assumption $A \geq C + \frac{3}{2}$. The other Equs.(A.12a-h) similarly imply that $A \geq C + \frac{3}{2}$ would give $t^{21}_{cd} = u^{12}_{ab} = u^{21}_{ab} = 0$, i.e., the $V_{\mu 12} = V_{\mu 21} = 0$.

Hence we conclude that both $C \geq A + \frac{3}{2}$ and $A \geq C + \frac{3}{2}$ are ruled out if $V_{\mu 12}, V_{\mu 21}$ are to be non-zero. In other words, given that they are integral or half-integral, $C - 1 \leq A \leq C + 1$. But we have already established that, for non-zero $V_{\mu 12}, V_{\mu 21}$ either $A$ is an integer and $C$ is a half-integer, or vice-versa. Consequently, the only remaining possible cases which avoid all the $V_{\mu}$ being zero are $A = C \pm \frac{1}{2}$.

The symmetry of Equs.(A.12a-h) under the exchanges $a \leftrightarrow c$, $b \leftrightarrow d$, $A \leftrightarrow C$, $B \leftrightarrow D$ means that $B = D \pm \frac{1}{2}$ is also required for non-zero $V_{\mu 12}, V_{\mu 21}$. There is no requirement for the sign in $B = D \pm \frac{1}{2}$ to be related to the sign in $A = C \pm \frac{1}{2}$. Pairs of blocks respecting the conditions $A = C \pm \frac{1}{2}$ and $B = D \pm \frac{1}{2}$ will be described as "compatible".

This almost completes the proof of Step 2 except that we need to show that when $A = C \pm \frac{1}{2}$ and $B = D \pm \frac{1}{2}$ a non-zero solution for the $V_{\mu 12}, V_{\mu 21}$ exists, and what it is. The fact that the four cases $A = C \pm \frac{1}{2}$ and $B = D \pm \frac{1}{2}$ (independent signs) do indeed produce non-trivial solutions can be seen as follows. Assume $C = A + \frac{1}{2}$. Since we require $a = c + \frac{3}{2}$ for (A.13) to hold we cannot assume $a = -A$ because this would imply an impossible c. Instead the smallest $a$ with which we can start the recursion is $a = -A + 1$. This gives,

$$r^{C}_{-A-\frac{1}{2}} t^{12}_{-A+1,b} = r^{A}_{-A} t^{12}_{-A,b} \tag{A.18}$$

But (4.4d) gives $r^{A}_{-A} = \sqrt{2A}$ and $r^{C}_{-A-\frac{1}{2}} = \sqrt{\left(C + A + \frac{1}{2}\right)\left(C - A + \frac{1}{2}\right)}$ which is non-zero for our assumed $C = A + \frac{1}{2}$. Hence we can take some non-zero $t^{12}_{-A,b}$ to start the recursion, in terms of which all the $t^{12}_{ab}$ are found for the assumed $b$. Similar recursion applies, of course, to the $b$ index, and this leads to all the $t^{12}_{ab}$ being found as factors times a single parameter, $t^{12}_{AB}$.

Similar remarks apply for $t^{21}_{cd}$, $u^{12}_{ab}$ and $u^{21}_{cd}$.

However, Shurtleff, Ref.[5], has also shown that the full set of recursion relations like (A.18) lead to $u^{12}_{-A,-B} = \pm t^{12}_{AB}$, etc., so there are in fact only two free parameters, $t^{12}_{AB}$ and $t^{21}_{CD}$. Moreover, in the notation, the $AB$ in $t^{12}_{AB}$ is implied by the 12, and the $CD$ in $t^{21}_{CD}$ is implied by the 21, so the two undetermined constants can be written more simply as $t_{12}$ and $t_{21}$. $V_{\mu 12}$ is proportional to $t_{12}$ whilst $V_{\mu 21}$ is proportional to $t_{21}$.

The derivation so far has focussed on $V_x$ and $V_y$ (or equivalently $V_\pm$) and their CRs with $J_i$. Bringing in the CRs with $K_i$ a similar derivation may be conducted for $V_t$ and $V_z$ (or the equivalent quantity $W_\pm = (V_z \pm V_t)/2$).

Shurtleff, Ref.[5], has thus derived the explicit expressions for all the $V_{\mu 12}$ and $V_{\mu 21}$ matrices, for all four combinations of cases $A = C \pm \frac{1}{2}$ and $B = D \pm \frac{1}{2}$. These are reproduced in Table A.1 but in a compact notation in which all four cases can be written as a single expression.

We have also adopted a different normalisation for the constants $t_{12}$ and $t_{21}$ which simplifies the expressions (and also simplifies their numerical values as will be seen later). The relationship between the $t_{12}$, $t_{21}$ used here and those of Ref.[5] are given in Table A.2 for completeness.

The matrices given in Table A.1 exclude the factors of $t_{12}$, $t_{21}$ and so are denoted as $U_{\mu 12}$ and $U_{\mu 21}$. The actual $V_{\mu 12}$ and $V_{\mu 21}$ matrices in any given case will therefore be given by,

$$V_{\mu 12} = t_{12} U_{\mu 12} \text{ and } V_{\mu 21} = t_{21} U_{\mu 21} \qquad \text{(A.19)}$$

And we have also defined $V_\pm = \left(V_x \pm iV_y\right)/2)$ and $W_\pm = (V_z \pm V_t)/2)$. Table A.1 gives expressions for the combinations $U_{\pm 12} = \left(U_{x12} \pm iU_{y12}\right)/2$ and $\tilde{U}_{\pm 12} = (U_{z12} \pm U_{t12})/2$, and similarly for $U_{\pm 21}, \tilde{U}_{\pm 21}$.

Finally, all the above analysis applies between any pair of (different) blocks, so 12 can be replaced throughout by arbitrary blocks $PQ$, and the signs in Table A.1, $S_A$ and $S_B$, align with,

$$A_P = A_Q + \frac{S_A}{2} \text{ and } B_P = B_Q + \frac{S_B}{2} \qquad \text{(A.20)}$$

**Table A.1: Explicit Expressions for the $U_{\mu 12}$ and $U_{\mu 21}$ Matrices** (based on Ref.[5])

| Quantity | Matrix Elements |
|---|---|
| $(U_{\pm 12})_{a_1 b_1, a_2 b_2}$ | $S_\pm \sqrt{(A_{12} \pm S_A a_{12})(B_{12} \pm S_B b_{12})} \cdot \delta_{a_1, a_2 \pm \frac{1}{2}} \delta_{b_1, b_2 \pm \frac{1}{2}}$ |
| $(U_{\pm 21})_{a_2 b_2, a_1 b_1}$ | $S_\pm \sqrt{(A_{12} \mp S_A a_{12})(B_{12} \mp S_B b_{12})} \cdot \delta_{a_1, a_2 \mp \frac{1}{2}} \delta_{b_1, b_2 \mp \frac{1}{2}}$ |
| $\left(\tilde{U}_{\pm 12}\right)_{a_1 b_1, a_2 b_2}$ | $-S_{AB}^{\pm} \sqrt{(A_{12} \pm S_A a_{12})(B_{12} \mp S_B b_{12})} \cdot \delta_{a_1, a_2 \pm \frac{1}{2}} \delta_{b_1, b_2 \mp \frac{1}{2}}$ |
| $\left(\tilde{U}_{\pm 21}\right)_{a_2 b_2, a_1 b_1}$ | $S_{AB}^{\pm} \sqrt{(A_{12} \mp S_A a_{12})(B_{12} \pm S_B b_{12})} \cdot \delta_{a_1, a_2 \mp \frac{1}{2}} \delta_{b_1, b_2 \pm \frac{1}{2}}$ |

Conventions/Definitions:

- Writing $A_1 = A, B_1 = B, A_2 = C, B_2 = D$ the signs $S_A, S_B$ are defined by $A = C + \frac{S_A}{2}$ and $B = D + \frac{S_B}{2}$. There are four cases which we may write $++, +-, -+, --$.
- $S_+ = +, S_- = -S_A S_B$.
- $S_{AB}^+ = S_B,\ S_{AB}^- = S_A$
- $A_{12}$ is the larger of $A$ or $C$ (or $A_1$ and $A_2$); $B_{12}$ is the larger of $B$ or $D$ (or $B_1$ and $B_2$).

- $a_{12}$ is $a_1$ or $a_2$ according to $A_{12}$ being $A_1$ or $A_2$. Similarly, $b_{12}$ is $b_1$ or $b_2$ according to $B_{12}$ being $B_1$ or $B_2$.
- $U_{\pm 12} = \left(U_{x12} \pm iU_{y12}\right)/2$ and $\widetilde{U}_{\pm 12} = (U_{z12} \pm U_{t12})/2$, similarly for $U_{\pm 21}, \widetilde{U}_{\pm 21}$.
- $J_\pm = J_x \pm iJ_y$ and $K_\pm = K_x \pm iK_y$

**Table A.2: Relationship between the $t_{12}$, $t_{21}$ used here and those of Ref.[5]**

| This work | Case ($S_A S_B$) | Equivalent in Ref.[5] |
|---|---|---|
| $t_{12}$ | $++$ | $\frac{t_{12}}{2\sqrt{A_1 B_1}}$ |
| | $--$ | $t_{12}$ |
| | $+-$ | $\frac{t_{12}}{\sqrt{2A_1}}$ |
| | $-+$ | $\frac{t_{12}}{\sqrt{2B_1}}$ |
| $t_{21}$ | $++$ | $t_{21}$ |
| | $--$ | $\frac{t_{21}}{2\sqrt{A_2 B_2}}$ |
| | $+-$ | $\frac{t_{21}}{\sqrt{2B_2}}$ |
| | $-+$ | $\frac{t_{21}}{\sqrt{2A_2}}$ |

### A.1 The Poincaré versus the dS/AdS Groups

It may now be checked directly that the explicit expressions of Table A.1, together with those for $J_i, K_i$, (4.4a-f), obey all the CRs of (2.1) and (2.2). These CRs are common between the Poincaré and dS/AdS groups and so the representations of Table A.1 and (4.4a-f) apply to both. Where the Poincaré and dS/AdS groups differ is in respect of CRs (2.3), the commutation relations of the $V_\mu$ amongst themselves. The only freedom left in order to constrain the representations of Table A.1 to obey these further CRs is the freedom to choose the parameters $t_{PQ}$. In the case of the Poincaré group, the corresponding CRs are $\left[V_\mu, V_\nu\right] = 0$, i.e., that the 4-momenta commute in flat Minkowski spacetime. This can be achieved only by requiring either $t_{PQ} = 0$ and, for compatible blocks, $t_{QP} \neq 0$, or vice-versa, so that the $V_\mu$ are neither Hermitian nor anti-Hermitian. In this work that case is ruled out because we are interested only in representations for which all generators are Hermitian or anti-Hermitian.

It follows that in the rest of this paper, provided $V_\mu$ of the form given by Table A.1 are used, the CRs of (2.1) and (2.2) will automatically be respected, and we need concentrate upon satisfying the CRs (2.3) only.

### A.2 Imposing the Hermitian/Anti-Hermitian Condition

Recall that the dSA requires $V_k^+ = V_k$ and $V_t^+ = -V_t$ and these require $V_+^+ = V_-$ as well as $W_+^+ = W_-$. But, considering the 12 and 21 blocks only, $V_+ = \begin{pmatrix} 0 & t_{12}U_{+12} \\ t_{21}U_{+21} & 0 \end{pmatrix}$ and $V_- = \begin{pmatrix} 0 & t_{12}U_{-12} \\ t_{21}U_{-21} & 0 \end{pmatrix}$. So, $V_+^+ = V_-$ requires $t_{12}^* U_{+12}^+ = t_{21} U_{-21}$. But for cases $++$ and $--$

Table A.1 tells us that $U^{+}_{+12} = -U_{-21}$ so for these cases $t^{*}_{12} = -t_{21}$ (e.g., if they are real they are of opposite sign; if imaginary they are of the same sign).

Alternatively, for cases $+-$ and $-+$ Table A.1 tells us that $U^{+}_{+12} = U_{-21}$ so for these cases $t^{*}_{12} = t_{21}$ (e.g., if they are real they are of the same sign; if imaginary they are of the opposite sign). These conditions are consistent with the Hermitian/anti-Hermitian requirements upon $W_{\pm}$ for the same reason. In summary,

For cases $++$ and $--$: $$t^{*}_{PQ} = -t_{QP} \tag{A.21a}$$

For cases $+-$ and $-+$: $$t^{*}_{PQ} = t_{QP} \tag{A.21b}$$

The simplicity of these requirements is an additional reason for adopting our present convention for the definitions of the $t_{PQ}$, differing from that of Ref.[5], see Table A.2.

## Appendix B: Construction of the Monotonic Irreps m/1 and 1/n

We have proved already in §6 that there must be two or more blocks in the backbone. With $J_i, K_i$ given in block-diagonal form, e.g., (A.1), we assert that $V_\mu$ consistent with CRs (2.3), and with backbones of the form m/1 or 1/n, can be written in the following manner… exemplifying this by the four block case for illustration,

$$\begin{pmatrix} 0 & V_{\mu12} & 0 & 0 \\ V_{\mu21} & 0 & V_{\mu23} & 0 \\ 0 & V_{\mu32} & 0 & V_{\mu34} \\ 0 & 0 & V_{\mu43} & 0 \end{pmatrix} \tag{B.1}$$

In other words that possible $V_\mu$ can be written with off-BD blocks which are BD-adjacent, i.e., within one block row and column of the BD. In the four-block case illustrated this means we can take $V_{\mu14} = V_{\mu41} = 0$. Note that the structure of the backbones for m/1 and 1/n, and the order in which the direct sums are written, ensures that adjacent blocks are compatible and hence can be taken as connected (unless it should turn out that the corresponding $V_{\mu PQ}$ is identically zero). It also follows, for types m/1 and 1/n, that no other off-BD blocks are non-zero.

The claim that $V_\mu$ of this form meet all the requirements to ensure a representation exists is to be proved in this section.

To prove the assertion it will suffice to show that $t_{12}, t_{23}, t_{34}, \ldots$ exist which respect the CRs (2.3), noting that $t_{21}, t_{32}, t_{43}, \ldots$ follow from (A.21a,b). This is by no means obvious because, on the face of it, the situation is hugely over-constrained. For example, in the four irrep case there are just 3 variables to be determined in order to satisfy 6 matrix equations, (2.3), where each matrix equation consists of 30 component equations (type 1/n) or 20 component equations (case m/1). Of course, these equations are not all independent.

Consider $[V_\alpha, V_\beta]$. We will illustrate the construction of the solution using the four block case, the extension to an arbitrary number of blocks will be clear. The on-BD blocks in $[V_\alpha, V_\beta]$ are,

$$[V_\alpha, V_\beta]_{11} = (V_{\alpha12}V_{\beta21} - V_{\beta12}V_{\alpha21}) \tag{B.2a}$$

$$[V_\alpha, V_\beta]_{22} = (V_{\alpha21}V_{\beta12} - V_{\beta21}V_{\alpha12}) + (V_{\alpha23}V_{\beta32} - V_{\beta23}V_{\alpha32}) \tag{B.2b}$$

$$[V_\alpha, V_\beta]_{33} = (V_{\alpha32}V_{\beta23} - V_{\beta32}V_{\alpha23}) + (V_{\alpha34}V_{\beta43} - V_{\beta34}V_{\alpha43}) \tag{B.2c}$$

$$[V_\alpha, V_\beta]_{44} = (V_{\alpha43}V_{\beta34} - V_{\beta43}V_{\alpha34}) \tag{B.2d}$$

As regards any larger number of blocks in the backbone, by assuming that only BD-adjacent blocks occur in the $V_\alpha$ there will only ever be one term like those in brackets in the first and last BD blocks, and only ever two such in the others. (In contrast, the inclusion of non-BD-adjacent blocks would bring in more such terms).

The (potentially) non-zero off-BD blocks are,

$$[V_\alpha, V_\beta]_{13} = (V_{\alpha12}V_{\beta23} - V_{\beta12}V_{\alpha23}) \tag{B.3a}$$

$$[V_\alpha, V_\beta]_{31} = (V_{\alpha32}V_{\beta21} - V_{\beta32}V_{\alpha21}) \tag{B.3b}$$

$$\left[V_\alpha, V_\beta\right]_{24} = \left(V_{\alpha 23}V_{\beta 34} - V_{\beta 23}V_{\alpha 34}\right) \tag{B.3c}$$

$$\left[V_\alpha, V_\beta\right]_{42} = \left(V_{\alpha 43}V_{\beta 32} - V_{\beta 43}V_{\alpha 32}\right) \tag{B.3d}$$

As regards any larger number of blocks in the backbone, by assuming that only BD-adjacent blocks occur in the $V_\alpha$ there will only ever be one term like those in brackets in all the off-BD blocks.

Because $J_i, K_i$ are in block-diagonal form, e.g., (A.1), the CRs (2.3) require all the non-BD blocks of $\left[V_\alpha, V_\beta\right]$ to be zero, so we shall require all of (B.3a-d) to be identically zero.

**B.1 Backbone Type 1/n**

This is the case $++, ++, ++, ....$ (It can equally be considered as $--, --, --, ...$, by reversing the order of the blocks). Evaluating $\left[V_x, V_y\right]_{13}$ using (B.3a) and the representations of Table A.1 shows it to be identically zero, and the same is found for all the $\left[V_\alpha, V_\beta\right]_{13}$. By symmetry the same is true for all the off-BD blocks, (B.3a-d). This confirms the CRs (2.3) are respected as regards the off-BD blocks.

Turning now to the on-BD blocks, we find,

$$i\left(\left[V_x, V_y\right]_{11}\right)_{a_1b_1,a_1'b_1'} = 4t_{12}t_{21}(A_1b_1 + B_1a_1)\delta_{a_1a_1'}\delta_{b_1b_1'} \tag{B.4}$$

This is required to be identically equal to $-(J_{z11})_{a_1b_1,a_1'b_1'} = -(a_1 + b_1)\delta_{a_1a_1'}\delta_{b_1b_1'}$, i.e.,

$$4t_{12}t_{21}(A_1b_1 + B_1a_1) \equiv -(a_1 + b_1) \tag{B.5}$$

We note that $A_1 > 0$ and $B_1 > 0$ because we have shown there must be at least two blocks in the backbone this case therefore requires $A_1 = B_1 \geq \frac{1}{2}$. Note that the minus sign in (B.5) is consistent with (A.21a) and we have,

$$t_{12} = -t_{21} = \frac{1}{2\sqrt{A}} = \sqrt{\frac{1}{2(N-1)}} \tag{B.6}$$

where $A_1 = B_1 = A = \frac{N-1}{2}$ and $N$ is the number of blocks in the backbone. It can similarly be checked by direct substitution of the representations of Table A.1 that,

$$\left[V_y, V_z\right]_{11} = iJ_x \qquad \text{and} \qquad \left[V_z, V_x\right]_{11} = iJ_y \tag{B.7}$$

But it is easier to simply note that both $J_i$ and $V_i$ transform as 3-vectors under rotations in $E^3$ and hence that (B.7) follows immediately from $\left[V_x, V_y\right]_{11} = iJ_z$.

Similarly, it can be checked directly from Table A.1 that the CRs $[V_t, V_i] = iK_i$ are obeyed, but it is easier to simply observe that one of these CRs was used in the derivation of the expressions of Table A.1 whilst the other two result from the 3-vector nature of $K_i$ and $V_i$.

We note in passing that (B.5) would have no solution in the case $A_1 \neq B_1$.

Turning now to the next on-BD block, $\left[V_x, V_y\right]_{22}$ evaluation using (B.2b) and Table A.1 gives,

$$i\left(\left[V_x, V_y\right]_{22}\right)_{a_2b_2,a_2'b_2'} = \tag{B.8}$$

$$\left[-4t_{12}t_{21}\left((A_2+1)b_2+(B_2+1)a_2\right)+4t_{23}t_{32}(A_2b_2+B_2a_2)\right]\delta_{a_2a_2'}\delta_{b_2b_2'}$$

This must be identically equal to $-(J_{z22})_{a_2b_2,a_2'b_2'}$=$-(a_2+b_2)\delta_{a_2a_2'}\delta_{b_2b_2'}$ hence,

$$-4t_{12}t_{21}\left((A_2+1)b_2+(B_2+1)a_2\right)+4t_{23}t_{32}(A_2b_2+B_2a_2)=-(a_2+b_2) \tag{B.9}$$

But $A_1=B_1=\frac{N-1}{2}$ and so $A_2=B_2=\frac{N-2}{2}$ and we have already found $t_{12}t_{21}$, (B.6), so we find that (B.9) is an identity provided we set,

$$t_{23}=-t_{32}=\sqrt{\frac{2N-1}{2(N-1)(N-2)}} \tag{B.10}$$

Block 22 of the other CRs of (2.3) will also be obeyed by the same rotational symmetry argument as before.

Turning now to the next on-BD block, $\left[V_x,V_y\right]_{33}$ it may be evaluated using (B.2c) and Table A.1 but can be written down immediately from (B.8), giving,

$$i\left(\left[V_x,V_y\right]_{33}\right)_{a_3b_3,a_3'b_3'}= \tag{B.11}$$

$$\left[-4t_{23}t_{32}\left((A_3+1)b_3+(B_3+1)a_3\right)+4t_{34}t_{43}(A_3b_3+B_3a_3)\right]\delta_{a_3a_3'}\delta_{b_3b_3'}$$

This must be identically equal to $-(J_{z33})_{a_3b_3,a_3'b_3'}$=$-(a_3+b_3)\delta_{a_3a_3'}\delta_{b_3b_3'}$. Using $A_2=B_2=\frac{N-2}{2}$ and $A_3=B_3=\frac{N-3}{2}$ together with (B.10) we find that it is indeed an identity provided we set,

$$t_{34}=-t_{43}=\sqrt{\frac{3(N-1)}{2(N-2)(N-3)}} \tag{B.12}$$

The solution for the general $t_{n,n+1}$ is now clear, for $2\le n\le N-1$, namely,

$$t_{n,n+1}=-t_{n+1,n}=\sqrt{\frac{(2N-n+1)n}{4(N-n)(N-n+1)}} \tag{B.13}$$

But finally we have the last on-BD block, $\left[V_x,V_y\right]_{NN}$, for which there is only one term as illustrated by (B.2d). Equating with $J_{zNN}$ we find,

$$-4t_{N,N-1}t_{\mathrm{N-1,N}}\left((A_N+1)b_{\mathrm{N}}+(B_{\mathrm{N}}+1)a_{\mathrm{N}}\right)\equiv-(a_N+b_N) \tag{B.14}$$

But because $t_{N,N-1}t_{\mathrm{N-1,N}}<0$ due to (A.21a) the LHS, if non-zero, would be positive when neither of $a_N$ or $b_N$ are negative and at least one is non-zero, whereas the RHS would then be negative. Hence (B.14) can only be an identity if both $a_N$ and $b_N$ are constrained to be zero, i.e., for $A_N=B_N=0$.

This establishes that, had we chosen the first block of the backbone, $(A,A)$, to have $A>\frac{N-1}{2}$, the final block of the backbone would not be $(0,0)$ and this would provide a conflict so that no solution would result. Similarly, if $A<\frac{N-1}{2}$ a conflict would again result. To see this, set $A=\frac{N'-1}{2}$ where $N'\le N-1$ in which case (B.13) becomes, for $n=N-1$,

$$-t_{n,n+1}t_{n+1,n}=\frac{(2N'-N+2)(N-1)}{4(N'-N+1)(N'-N+2)}$$

But this is undefined due to the zero in the denominator when $N' = N - 1$. For smaller $N'$ this same problem will occur for smaller $n$, i.e., $n < N - 1$. Hence we must have $A = \frac{N-1}{2}$, no other option is consistent.

The on-BD blocks of the other CRs (2.3) can be established for backbone 1/n in similar fashion, though, for the same argument noted above, there is no need to do so as this follows from the derivation of the representations of Table A.1 together with the 3-vector nature of $J_i, K_i$ and $V_i$. Table B.1 gives explicitly the values of the coefficients $t_{ij}$ for the first ten irreps.

**Summary of §B.1**

It has been established by explicit construction of the values of the $t_{n,n+1} = -t_{n+1,n}$, (B.13), that a representation of the dSA exists on the backbone of 1/n. The solution was based on the assumption that the representations of the $V_\mu$ have only BD-adjacent blocks, which corresponds to the connectivity of Figure 4b. As a by-product we have also proved that there is no solution on a backbone of the following type,

$(A, B) \oplus \left(A - \frac{1}{2}, B - \frac{1}{2}\right) \oplus (A - 1, B - 1) \oplus \dots \oplus \left(A - \frac{N-1}{2}, B - \frac{N-1}{2}\right)$

unless it reduces to 1/n, i.e., that $A = B = \frac{N-1}{2}$ so that the last block is $(0,0)$.

### B.2 Backbone Type m/1

In this section we seek to carry out a similar derivation as in §B.1 for the backbone m/1. This is the case $+-, +-, +-, \dots$ (which, reading the bocks in reverse order is equivalent to $-+, -+, -+, \dots$).

Evaluating $[V_x, V_y]_{13}$ using (B.3a) and the representations of Table A.1 shows it to be identically zero, and the same is found for all the $[V_\alpha, V_\beta]_{13}$. To show this it is not necessary to assume that $B_1 = 0$, or any relation between $A_1$ and $B_1$, and so the same must be true for all the off-BD blocks, (B.3a-d). This confirms the CRs (2.3) are respected as regards the off-BD blocks.

Turning now to the on-BD blocks, we find, using (B.2a) and Table A.1,

$$i\left([V_x, V_y]_{11}\right)_{a_1 b_1, a_1' b_1'} = 4t_{12}t_{21}(A_1 b_1 - (B_1 + 1)a_1)\delta_{a_1 a_1'}\delta_{b_1 b_1'} \qquad \text{(B.15)}$$

and this is required to be identically equal to $-(J_{z11})_{a_1 b_1, a_1' b_1'} = -(a_1 + b_1)\delta_{a_1 a_1'}\delta_{b_1 b_1'}$, i.e.,

$$4t_{12}t_{21}(A_1 b_1 - (B_1 + 1)a_1) = -(a_1 + b_1) \qquad \text{(B.16)}$$

But $A_1 = A_2 + \frac{1}{2} > 0$ and so if $B_1 > 0$ (B.16) could not be identically true. Hence a backbone which starts $(A_1, B_1) \oplus \left(A_1 - \frac{1}{2}, B_1 + \frac{1}{2}\right)$ must have $B_1 = 0$. (B.16) then gives, recalling that (A.21b) requires $t_{21} = t_{12}$,

$$t_{21} = t_{12} = \frac{1}{2} \qquad \text{(B.17)}$$

For the next on-BD block, we find, using (B.2b) and Table A.1, (B.18)

$$i\left(\left[V_x,V_y\right]_{22}\right)_{a_2b_2,a_2'b_2'} = \left[4t_{23}t_{32}(A_2b_2-(B_2+1)a_2)-\left((A_2+1)b_2-B_2a_2\right)\right]\delta_{a_2a_2'}\delta_{b_2b_2'}$$

and this is required to be identically equal to $-(J_{z22})_{a_2b_2,a_2'b_2'} = -(a_2+b_2)\delta_{a_2a_2'}\delta_{b_2b_2'}$, i.e.,

$$4t_{23}t_{32}(A_2b_2-(B_2+1)a_2)-\left((A_2+1)b_2-B_2a_2\right) = -(a_2+b_2) \tag{B.19}$$

We find that (B.19) is an identity iff,

$$t_{32} = t_{23} = \frac{1}{2} \tag{B.20}$$

For the next on-BD block, we find, in similar fashion, (B.21)

$$i\left(\left[V_x,V_y\right]_{33}\right)_{a_3b_3,a_3'b_3'} = \left[4t_{34}t_{43}(A_3b_3-(B_3+1)a_3)-\left((A_3+1)b_3-B_3a_3\right)\right]\delta_{a_3a_3'}\delta_{b_3b_3'}$$

and this is required to be identically equal to $-(J_{z22})_{a_3b_3,a_3'b_3'} = -(a_2+b_2)\delta_{a_3a_3'}\delta_{b_3b_3'}$, i.e.,

$$4t_{34}t_{43}(A_3b_3-(B_3+1)a_3)-\left((A_3+1)b_3-B_3a_3\right) = -(a_3+b_3) \tag{B.22}$$

Again this is an identity iff,

$$t_{43} = t_{34} = \frac{1}{2} \tag{B.23}$$

For the case of a backbone with $N$ blocks of type $(A,0)\oplus\left(A-\frac{1}{2},\frac{1}{2}\right)\oplus(A-1,1)\oplus$ ... it is clear that from block $\left[V_x,V_y\right]_{nn}$ with $n \le N-1$ we will find,

$$t_{n+1,n} = t_{n,n+1} = \frac{1}{2} \tag{B.24}$$

That leaves only block $i\left[V_x,V_y\right]_{NN}$ which differs by containing only one term, see (B.2d). Equating with $-J_{zNN}$ requires,

$$-4t_{N-1,N}t_{N,N-1}\left((A_N+1)b_N-B_Na_N\right) = -(a_N+b_N) \tag{B.25}$$

But as $B_N > 0$ if $A_N > 0$ then (B.25) cannot be an identity. Hence it must be that $A_N = 0$, i.e., the backbone must terminate at $A_N = 0$. (B.25) is then an identity as long as we have,

$$t_{N-1,N} = t_{N,N-1} = \frac{1}{2} \tag{B.26}$$

which is consistent with (B.24), hence completing the proof that this representation exists, noting that the backbone terminates at $A_N = 0$ also means $A_1 = \frac{N-1}{2}$.

**Summary of §B.2**

It has been established by explicit demonstration that $t_{n+1,n} = t_{n,n+1} = \frac{1}{2}$ provides a representation of the dSA on the backbone of m/1. The solution was based on the assumption that the representations of the $V_\mu$ have only BD-adjacent blocks, consistent with the connectivity of Figure 4a. As a by-product we have also proved that there is no solution on a backbone of the following type,

$$(A,B) \oplus \left(A-\tfrac{1}{2}, B+\tfrac{1}{2}\right) \oplus (A-1, B+1) \oplus \ldots \oplus \left(A-\tfrac{N-1}{2}, B+\tfrac{N-1}{2}\right)$$

unless it reduces to type m/1, i.e., that $A = \frac{N-1}{2}$ and $B = 0$ so that the first block is $(A, 0)$ and the last block is $(0, A)$.

We shall refer to the representations 1/n and m/1, whose existence has been proved in this Appendix, as irreps as this is demonstrated in §11 by evaluation of the Casimir matrices.

**Table B.1: The Coefficients $t_{ij}$ which Uniquely Specify the First Ten Irreps.** Only the irreps of Type 1/n are listed below and these have $t_{ji} = -t_{ij}$. For irreps of Type m/1, $t_{ij} = \frac{1}{2}$ for all $t_{n,n+1}$ for $n \in [1, N-1]$.

| Rep | $t_{12}$ | $t_{23}$ | $t_{34}$ | $t_{45}$ | $t_{56}$ |
|---|---|---|---|---|---|
| 1/2 | $\frac{1}{\sqrt{2}}$ | - | - | - | - |
| 1/3 | $\frac{1}{2}$ | $\frac{\sqrt{5}}{2}$ | - | - | - |
| 1/4 | $\frac{1}{\sqrt{6}}$ | $\sqrt{\frac{7}{12}}$ | $\frac{3}{2}$ | - | - |
| 1/5 | $\frac{1}{\sqrt{8}}$ | $\sqrt{\frac{3}{8}}$ | 1 | $\sqrt{\frac{7}{2}}$ | - |
| 1/6 | $\frac{1}{\sqrt{10}}$ | $\sqrt{\frac{11}{40}}$ | $\sqrt{\frac{5}{8}}$ | $\sqrt{\frac{3}{2}}$ | $\sqrt{5}$ |

## Appendix C: The Characterisation of Reducible Representations and the Significance of Duplicate Blocks

Here we make some elementary observations on what characterises reducible representations on the $A, B$ plots used in the main text to characterise the irreps. We also observe that if a representation contains a duplicate block, then it must be reducible, i.e., the irreps have no duplicate blocks. (A backbone contains a "duplicate" if there are two or more blocks on the backbone with the same $(A, B)$).

A backbone which consists of the backbones of two or more different irreps but with non-intersecting or overlapping characterising diagrams clearly provides a representation which is the direct sum of the corresponding irreps, i.e., it is reducible, being already in reduced form. Figures C.1 is an example.

**Figure C.1: An Example Backbone for a Reducible Representation.** *Only shaded boxes are present, others shown below are merely unoccupied place holders. Blocks which are connected are shown linked. Hence this is the backbone for the reducible representation based on the direct sum of the irrep based on* $(2,0) \oplus \left(\frac{3}{2}, \frac{1}{2}\right) \oplus (1,1) \oplus \left(\frac{1}{2}, \frac{3}{2}\right) \oplus (0,2)$ *and that based on* $(0,0) \oplus \left(\frac{1}{2}, \frac{1}{2}\right)$*. These are bosonic cases.*

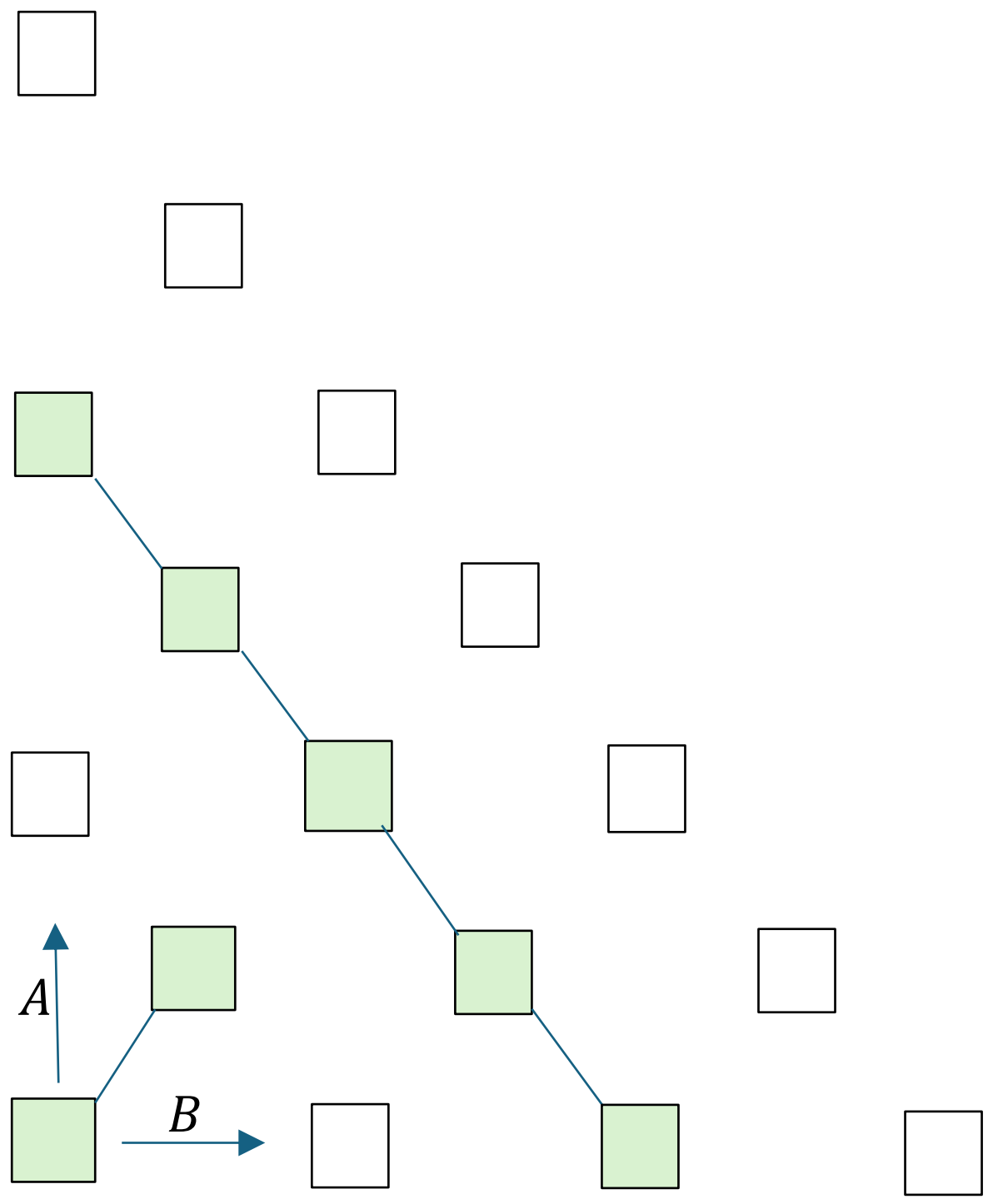

But what happens if trajectories intersect? For example consider that shown in Figure C.2. This fails to provide either an irrep or a reducible representation because it contains a unique non-monotonic 3-sequence, 123, and another, 124. Hence by §8 it is not a representation.

However, we can modify Figure C.2 by introducing a duplicate block at position 2, as shown in Figure C.3, the duplicate block being labelled 4.

**Figure C.2: An Example of Intersecting Backbone Trajectories.** *Only shaded boxes are present, others shown below are merely unoccupied place holders. Blocks which are connected are shown linked. Numbers are block labels. There is no representation on the backbone illustrated here.*

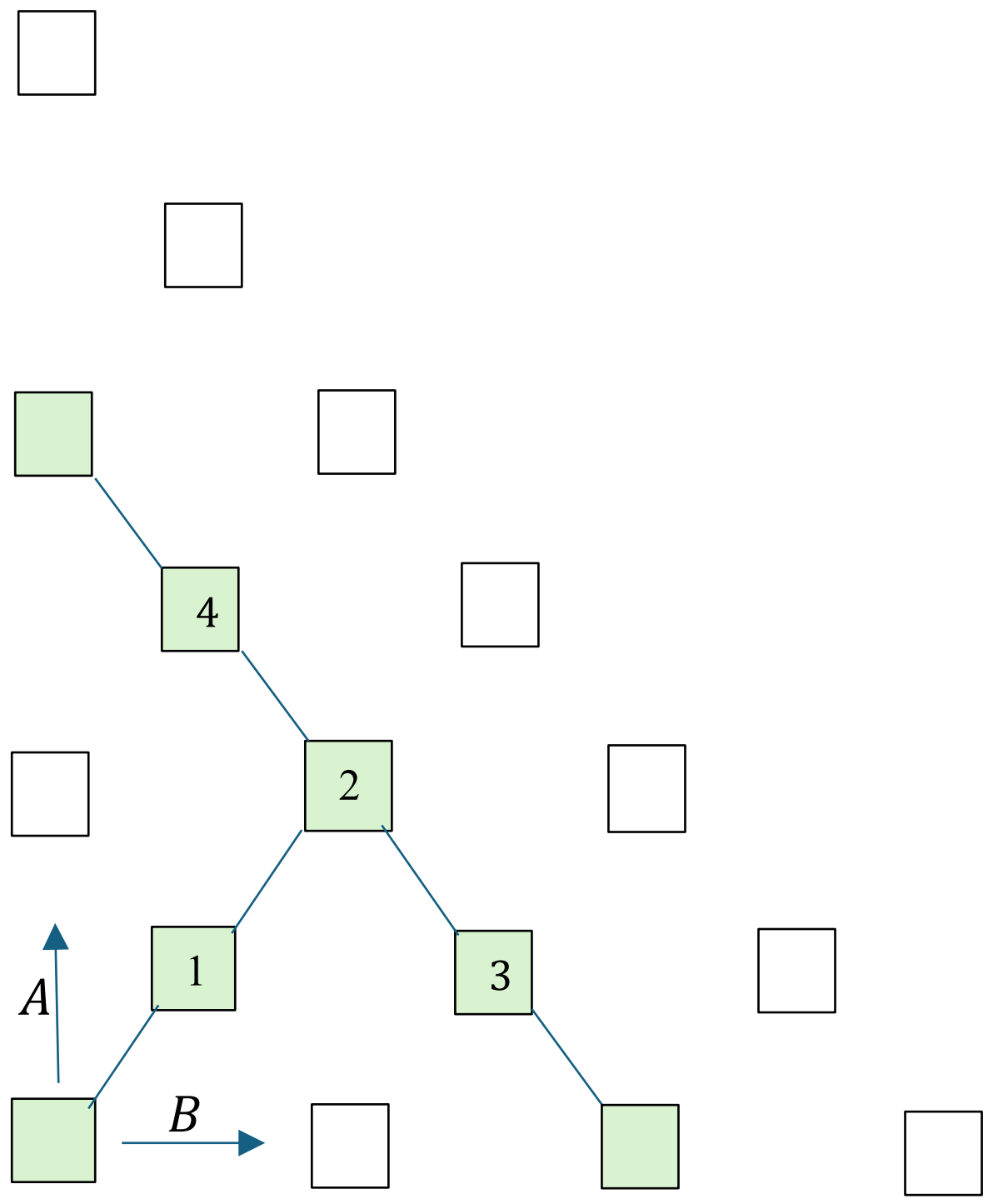

**Figure C.3: A Modified Version of the Example Case of Figure B.** *There is now a duplicate block, number 4, at the same* $(A, B)$ *as block 2 Numbers are block labels. A reducible representation results.*

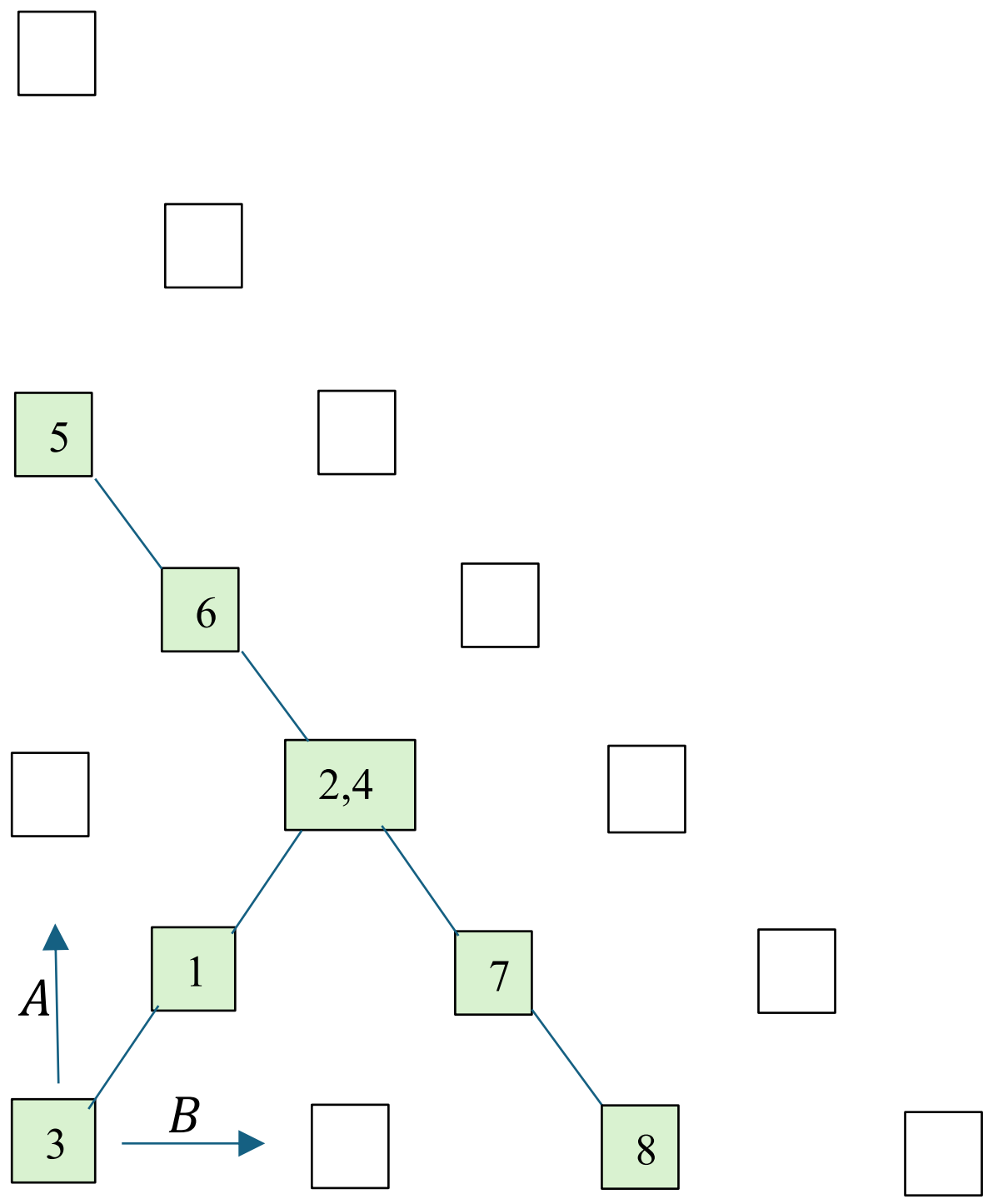

In Figure C.3 we take block 1 to be connected to block 4 but not to block 2, whilst blocks 6 and 7 are connected to block 2 but not to block 4. There is now no instance of a unique non-monotonic 3-sequence, we merely have a representation which provides both the irrep on backbone $(0,0) \oplus \left(\frac{1}{2},\frac{1}{2}\right) \oplus (1,1)$ and that on $(2,0) \oplus \left(\frac{3}{2},\frac{1}{2}\right) \oplus (1,1) \oplus \left(\frac{1}{2},\frac{3}{2}\right) \oplus (0,2)$, i.e., Figure C.3 provides a valid representation but it is reducible.

These informal observations illustrate how duplicates in a representation will indicate reducibility, whilst intersecting of overlapping diagrams without duplicates indicate no such representation.

## Appendix D: Uniqueness: Any Backbone has at Most One Representation (up to Equivalence)

It follows from the construction of the irreps in the main text that there is only one connectivity for a given backbone which provides an irrep. However, for completeness, we present here a more formal demonstration that any proposed backbone, if it produces a representation at all, produces only one, up to equivalence. In other words, if a given backbone has two differently realised representations for the $V_\mu$ then the two representations are equivalent (i.e., a unitary transformation exists to change one into the other).

The essence of the proof is simple and relies on the fact that irreps of the dSA can be constructed on a backbone of direct sums of the irreps of the algebra of the $SO(4)$ group, comprised of the generators $J_i, V_i$, instead of using, as we have here, backbones of direct sums of the irreps of the HLA comprised of the generators $J_i, K_i$. The key observation is that both algebras, up to complexification, are locally isomorphic to $\mathfrak{su}(2) \oplus \mathfrak{su}(2)$. Consequently, a given representation expressed in terms of an HLA backbone can be unitarily transformed into a representation based on another $\mathfrak{su}(2) \oplus \mathfrak{su}(2)$ backbone, that for the $J_i, V_i$.

If the representation is reducible, and hence expressible in block diagonal form over irreps, a unitary transformation between the two irrep bases can be chosen so as to preserve the block diagonal structure, i.e., the number of irreps, and the dimensionality of each one, is independent of the irrep basis used. Suppose the backbone in the HLA basis is given by $N$ blocks $(A_i, B_i)$ and that in the $SO(4)$ basis by $N$ blocks $(A_i', B_i')$. We thus have $(2A_i + 1)(2B_i + 1) = (2A_i' + 1)(2B_i' + 1)$ for all $i$. But the $SO(3)$ subgroup generated by the $J_i$ is the same regardless of the basis, and so each block $(A_i', B_i')$ must provide the same $J$ eigenvalues as $(A_i, B_i)$, hence $A_i + B_i = A_i' + B_i'$ and $|A_i - B_i| = |A_i' - B_i'|$. This is sufficient to imply that, for every $i$ independently, either $A_i' = A_i$ and $B_i' = B_i$ or alternatively that $A_i' = B_i$ and $B_i' = A_i$.

The case that $A_i' = B_i$ and $B_i' = A_i$ corresponds, via (2.4b), to the same $J_i$ but a $K_i$ which has changed sign. All the CRs, (2.1-3), are invariant under this change of sign, and hence the $i^{th}$ irrep is equivalent whether $A_i' = A_i$ and $B_i' = B_i$ or $A_i' = B_i$ and $B_i' = A_i$. If $A_i' = B_i$ and $B_i' = A_i$ applied for all $i$ then the overall representation would be equivalent to that using $A_i' = A_i$ and $B_i' = B_i$. However, if $A_i' = B_i$ and $B_i' = A_i$ applies for some but not all the $N$ irreps, and hence $A_i' = A_i$ and $B_i' = B_i$ applies for the rest, then overall this is a representation which is not equivalent to one based on $A_i' = A_i$ and $B_i' = B_i$ across all irreps. However, there is a further constraint to take into account, namely that the connectedness between blocks must be preserved. For example, if blocks 1 and 2 are connected then $A_1 = A_2 \pm \frac{1}{2}$ and $B_1 = B_2 \pm \frac{1}{2}$. But if $A_1' = B_1$ and $B_1' = A_1$ whilst $A_2' = A_2$ and $B_2' = B_2$ we would have only that $A_1' = B_2' \pm \frac{1}{2}$ and $B_1' = A_2' \pm \frac{1}{2}$ and connectedness in the alternative basis, which requires $A_1' = A_2' \pm \frac{1}{2}$ and $B_1' = B_2' \pm \frac{1}{2}$, is not automatically preserved. Careful consideration of the possible combinations of signs shows that connectedness is preserved only in two cases,

- $A_1 = B_1 = A_1' = B_1'$ and $A_2 = B_2 = A_2' = B_2'$, in which case the permutation $A \leftrightarrow B$ makes no difference, or,

- $(A_1', B_1') = (A_1 - 1, B_1 + 1)$ or $(A_1', B_1') = (A_1 + 1, B_1 - 1)$, but this case fails to preserve the minimum $J$ eigenvalue, as we get $|A_1{}' - B_1{}'| = |A_1 - B_1 \pm 2|$, and so can be ruled out.

It is concluded that the two backbone structures, $(A_i, B_i)$ based on the HLA, and $(A_i{}', B_i{}')$ based on SO(4), are identical. Thus any two realisations of the $V_\mu$ on the same backbone in the HLA basis have a unitary transformation to the same SO(4) basis backbone, and hence must have a unitary transformation which connects them, i.e., they are equivalent. In short, the representation on any given backbone, if it exists, is unique, up to equivalence.